\documentclass[12pt,twoside,a4paper]{article}
\usepackage[english]{babel}
\usepackage[utf8]{inputenc}
\usepackage{amsmath}
\usepackage{amsfonts}
\usepackage[pagebackref]{hyperref}
\usepackage{amsthm}
\usepackage{amssymb}
\usepackage{thmtools}
\usepackage{bm}
\usepackage{tikz-cd}
\usepackage{xr}
\usepackage[nameinlink]{cleveref}
\usepackage{url}

\usepackage{enumitem}
\usepackage{authblk}
\usepackage[title]{appendix}
\usepackage{chngcntr}
\usepackage{apptools}
\usepackage{mathtools}
\usepackage{appendix}

\usepackage{bookmark}
\usepackage{appendix}
\definecolor{crimsonglory}{rgb}{0.75, 0.0, 0.2}
\definecolor{darkpowderblue}{rgb}{0.0, 0.2, 0.6}
\hypersetup{
	colorlinks   = true, 
	urlcolor     = darkpowderblue, 
	linkcolor    = darkpowderblue, 
	citecolor   = crimsonglory 
}

\theoremstyle{plain}
\newtheorem{theorem}{Theorem}[section] 

\newtheorem{definition}[theorem]{Definition} 
\newtheorem{prop}[theorem]{Proposition}
\newtheorem{cor}[theorem]{Corollary}
\newtheorem{lemma}[theorem]{Lemma}

\newtheorem{remark}[theorem]{Remark}
\newtheorem{remarks}[theorem]{Remarks}

\newtheorem{assumption}[theorem]{Assumption}

\DeclareMathOperator{\End}{End}
\DeclareMathOperator{\GL}{GL}
\DeclareMathOperator{\SL}{SL}

\DeclareMathOperator{\spec}{Spec}

\DeclareMathOperator{\aut}{Aut}

\DeclareMathOperator{\disc}{disc}

\newcommand{\Sum}[2]{\displaystyle\sum_{#1}^{#2}}

\newcommand{\Z}{\mathbb{Z}}
\newcommand{\N}{\mathbb{N}}
\newcommand{\C}{\mathbb{C}}
\newcommand{\Q}{\mathbb{Q}}
\newcommand{\B}{\mathbb{B}}
\newcommand{\A}{\mathbb{A}}
\newcommand{\R}{\mathbb{R}}

\newcommand{\G}{\mathbb{G}}

\newcommand{\CE}{\mathcal{E}}
\newcommand{\CX}{\mathcal{X}}
\newcommand{\CG}{\mathcal{G}}
\newcommand{\ld}{,\ldots,}
\newcommand{\CP}{\mathcal{P}}
\newcommand{\CO}{\mathcal{O}}

\usepackage[OT2,T1]{fontenc}
\DeclareSymbolFont{cyrletters}{OT2}{wncyr}{m}{n}
\DeclareMathSymbol{\Sha}{\mathalpha}{cyrletters}{"58}

\usepackage{tocloft}

\makeatletter
\newcommand{\address}[1]{\gdef\@address{#1}}
\newcommand{\email}[1]{\gdef\@email{\url{#1}}}
\newcommand{\@endstuff}{\par\vspace{\baselineskip}\noindent\small
	\begin{tabular}{@{}l}\scshape\@address\\\textit{E-mail address:} \@email\end{tabular}}
\AtEndDocument{\@endstuff}
\makeatother

\title{Effective Brauer-Siegel on some curves in $Y(1)^n$}
\author{Georgios Papas}
\date{}
\address{Institute for Advanced Study\\
	1 Einstein Drive\\
	Princeton, N.J. 08540\\
	U.S.A.}
\email{gpapas@ias.edu}

\makeatletter
\newcommand{\subjclass}[2][1991]{%
	\let\@oldtitle\@title%
	\gdef\@title{\@oldtitle\footnotetext{#1 \emph{Mathematics subject classification.} #2}}%
}

\begin{document}
	\maketitle
	
	\begin{abstract}We establish an effective version of Siegel's lower bounds for class numbers of imaginary quadratic fields in certain curves in $Y(1)^n$. Our proof goes through the G-functions method of Yves Andr\'e. Following recent results of G. Binyamini, these lead to effective Andr\'e-Oort statements for the curves in question.

 \end{abstract}

	
	\section{Introduction}

\subsection{Siegel's lower bounds for class numbers}

In $1935$ Siegel established the following lower bounds on the class numbers of imaginary quadratic fields.

\begin{theorem}[Siegel, \cite{siegelog}]\label{siegelogthm}Given $\epsilon>0$ there exists a constant $c(\epsilon)>0$ such that 
	\begin{equation}\label{eq:siegeleq}
		h(K)\geq c(\epsilon) |\disc(K)|^{\frac{1}{2}-\epsilon},
	\end{equation}for all imaginary quadratic fields $K$, where $h(K)$ denotes the class number of the field $K$.
\end{theorem}

We note that it is known, see for example the discussion in $\S 2.2$ of \cite{galeffectiveandreoort}, that \eqref{eq:siegeleq} is known to hold for all orders in imaginary quadratic fields. In other words, given $\epsilon>0$, a constant $c(\epsilon)>0$ exists such that\begin{equation}
	h(\CO)\geq c(\epsilon) |\disc(\CO)|^{\frac{1}{2}-\epsilon}\end{equation}
	for all orders $\CO$ in any imaginary quadratic field $K$.

Due to the dependence of the proof on Siegel zeroes, the constants $c(\epsilon)$ that appear in \Cref{siegelogthm} are famously ineffective. Subsequent attempts to give effective proofs of \Cref{siegelogthm} have only been successful when replacing $|\disc(K)|^{\frac{1}{2}-\epsilon}$ in \eqref{eq:siegeleq} by the smaller quantity $(\log|\disc(K)|)^{1-\epsilon}$, due to work of Goldfeld-Gross-Zagier, see \cite{goldfeldsurvey}.

The closest to an effective version of \Cref{siegelogthm} is the following result of T. Tatuzawa:
\begin{theorem}[Tatuzawa, \cite{tatuzawa}]Given $\epsilon>0$ there exist an effective constant $c(\epsilon)>0$ and an imaginary quadratic field $K_{\epsilon}$ such that \eqref{eq:siegeleq} holds for all imaginary quadratic fields with $K\neq K_{\epsilon}$.
\end{theorem}
We note that the ``exceptional'' quadratic field $K_{\epsilon}$ in Tatuzawa's result, at least for small enough $\epsilon>0$, is likely to not exist at all. It is worth, perhaps, pointing out here that were we to assume the validity of the generalized Riemann hypothesis for imaginary quadratic fields, Siegel's result would be effective.

We note that Siegel's motivation for \Cref{siegelogthm} stemmed from an attempt to answer Gauss's class number problem, for a beautiful summary of the history of this problem and the connection to \Cref{siegelogthm} see \cite{goldfeldsurvey}. Also, for a simple proof of \Cref{siegelogthm} we point the interested reader to \cite{goldfeldsimple}.

\subsubsection{Connection to CM elliptic curves}\label{section:connectingtocmellcurves}

\Cref{siegelogthm} can be naturally connected to the theory of elliptic curves with complex multiplication. 

Consider $\tau\in \mathbb{H}$ in the complex upper half plane which is imaginary quadratic and let $E_{j(\tau)}$ be the elliptic curve with $j$ invariant $j(\tau)$. It is then classical that $E_{j(\tau)}$ admits complex multiplication by an order of the imaginary quadratic field $\Q(\tau)$, or in other words $\End(E_{j(\tau)})$ is an order in the imaginary quadratic field $\Q(\tau)$. 

Let us set, for $\tau$ as above, $D_\tau$ to be the discriminant of the primitive quadratic polynomial with coefficients in $\Z$ that has $\tau$ as a root. Also classical in this case, see for example \cite{cmseminar}, is the fact that the class number of the order $\End(E_{j(\tau)})$ is equal to $[\Q(\tau,j(\tau)):\Q(\tau)]=[\Q(j(\tau)):\Q]$. We can thus rephrase \Cref{siegelogthm} in the following format:
\begin{theorem}[Siegel's Theorem V. $2$]\label{siegelalternate} Given $\epsilon>0$ there exists a constant $c(\epsilon)>0$ such that for all imaginary quadratic $\tau \in \mathbb{H}$ we have 
		\begin{equation}\label{eq:siegeleqalt}
		[\Q(j(\tau)):\Q] \geq c(\epsilon) | \disc(\End(E_{j(\tau)})) |^{\frac{1}{2}-\epsilon}.
	\end{equation}	
\end{theorem}

\subsection{Our main result}

Throughout our exposition, $Y(1)$ will denote the level-one modular curve and $X(1)$ its compactification. These curves may be identified with the affine line $\mathbb{A}^1$, and projective line $\mathbb{P}^1$ respectively, via the isomorphism yielded by Klein's $j$-invariant. The main result of this text can be seen as an effective version of \Cref{siegelalternate} for tuples in certain $1$-parameter families in $Y(1)^n$. 

By a ``\textbf{CM point in }$Y(1)^n$'' we will mean a tuple $(\zeta_1\ld\zeta_n)$ such that, for all $i$, $\zeta_i$ is the $j$-invariant of an elliptic curve with complex multiplication. In this terminology, the reformulation \Cref{siegelalternate} is really a translation of Siegel's result to a statement about the size of the fields of definition of CM points in $Y(1)$.

We start with a bit of notation. Let $C\subset Y(1)^n$, where $n\geq 2$, be a smooth irreducible curve defined over $\bar{\Q}$ and let $\bar{C}$ be its Zariski closure in $X(1)^n$. We also let $s_0\in \bar{C}(\bar{\Q})\backslash Y(1)^n$ be a fixed point in the boundary $X(1)^n\backslash Y(1)^n$. 

\begin{definition}\label{defcoord}Let $C$, $s_0$ be as above and let $\pi_i:X(1)^n\rightarrow X(1)$ denote the coordinate projections. 
	
	The coordinate $i$ will be called \textbf{smooth for} $C$ \textbf{with respect to }$s_0$ if $\pi_i(s_0)\in Y(1)$. A smooth coordinate $i$ for the curve $C$ will be called a \textbf{CM coordinate} for $C$ \textbf{with respect to }$s_0$ if in addition $\pi(s_0)$ is a CM point in $Y(1)$. Finally, the coordinate $i$ will be called \textbf{ singular for }$C$ \textbf{with respect to }$s_0$ if it is not smooth, i.e. if $\pi_i(s_0)=\infty$.\end{definition}

Our main result is then the following:
\begin{theorem}\label{maintheorem}Let $C\subset Y(1)^n$ be a smooth irreducible curve defined over $\bar{\Q}$ that is not contained in a proper special subvariety of $Y(1)^n$. Let $\bar{C}$ be the compactification of $C$ in $X(1)^n$ and fix some $s_0\in \bar{C}(\bar{\Q})\backslash Y(1)^n$. Assume that with respect to $s_0$ there exists at least one CM coordinate for $C$ or at least two singular coordinates for $C$.
	
	Then there exist effectively computable positive constants $C_1$ and $C_2$, depending on $C$ and $s_0$, such that for all CM points $s=(s_1\ld s_n)\in C(\bar{\Q})$ we have \begin{center}
		$[\Q(s):\Q]\geq C_1 \max\{ |\disc(\End(E_{s_k}))|:1\leq k\leq n\}^{C_2}  $.
	\end{center}
\end{theorem}
\begin{remark}1. Upon identifying $Y(1)^n$ with the affine $n$-dimensional space $\spec(\Q[x_1\ld x_n])$ we may give a rough description of the special subvarieties of $Y(1)^n$. 
	
	For $N\in\N$ we denote by $\Phi_N(X,Y)\in\Z[X,Y]$ the $N$-th modular polynomial. Special subvarieties are the irreducible components of subvarieties of $Y(1)^n$ defined by systems of equations of the form \begin{enumerate}
	\item $\Phi_N(x_i,x_j)=0$, and
	\item $x_j=\zeta$,
\end{enumerate}where $\zeta$ is a singular modulus, i.e. $\zeta$ is the $j$-invariant of a CM elliptic curve. For a more precise description, we point the interested reader to Definition $1.3$ in the introduction of \cite{pilaao}.\\

2. The above result is more closely in spirit to an earlier result of E. Landau, see \cite{landau}. In more detail, Landau shows the existence of an ineffective constant $c$ for which 
\begin{equation}
			h(K)\geq c |\disc(K)|^{\frac{1}{8}},
\end{equation}for all imaginary quadratic fields $K$.
\end{remark}

The proof of our result follows Yves Andr\'e's ``G-functions method'', which we review in this setting in \Cref{section:background}. In particular, \Cref{maintheorem} will follow from height bounds for CM points on curves of the form considered in \Cref{maintheorem}. For some crude estimates on the constants that appear in \Cref{maintheorem} we point the interested reader to \Cref{remarkonconstants}.	
	\subsubsection{Connection to Unlikely Intersections}
	The author's motivation in pursuing these results stems from the field of unlikely intersections. In particular, our motivation comes from Pila's proof of the Andr\'e-Oort conjecture in $Y(1)^n$, see \cite{pilaao}. In his proof, Pila employs the so-called ``Pila-Zannier strategy''. The key arithmetic input in Pila's proof is Siegel's \Cref{siegelogthm}, or rather the alternate formulation \Cref{siegelalternate}, which shows that the Galois orbits of CM points are large enough for the strategy to work.
	
	Pila's usage of Siegel's lower bounds on class numbers \Cref{siegelalternate} is one of the reasons that his proof of Andr\'e-Oort is not effective. The same lower bounds were used by Andr\'e in \cite{andrefirstao} in establishing the Andr\'e-Oort conjecture for $Y(1)^2$. Effective proofs of this result of Andr\'e were later given by K\"uhne \cite{kuhneao} and Bilu-Masser-Zannier \cite{bilumasserzannier}, without using the ineffective lower bounds of Siegel.
	
	In $13.3$ of \cite{pilaao}, J. Pila notes that effectivity in the Andr\'e-Oort setting in $Y(1)^n$ would follow from the, conjectural at the time, effective o-minimality of the structure $\R_j$, where $j$ here denotes Klein's $j$-function, and effective lower bounds for Galois orbits. The effectivity of the structure $\R_j$ was recently established by G. Binyamini, see $\S 10.4$ in \cite{binyaminilog}. In the particular case of $Y(1)^n$ we note that Gal Binyamini's earlier papers \cite{galeffectiveandreoort,binyaminifoliations} circumvent the issue of the effectivity of the structure $\R_j$ and establish an effective version of Andr\'e-Oort modulo the ineffectivity of Siegel's result \Cref{siegelogthm}. All in all, as a result of our \Cref{maintheorem} one gets:
\begin{cor}\label{corollaryintro}Let $C\subset Y(1)^n$ be a smooth irreducible curve defined over $\bar{\Q}$ that is not contained in a proper special subvariety, and let $\bar{C}$ be the compactification of $C$ in $X(1)^n$.
	
	Assume that there exists $s_0\in\bar{C}(\bar{\Q})\backslash Y(1)^n$ with respect to which $C$ has at least two singular coordinates, or $C$ has at least one CM coordinate. Then the number of CM points in $C(\bar{\Q})$ can be effectively bounded.
\end{cor}

\Cref{corollaryintro} is a special case of the aforementioned results of L. K\"uhne, see \cite{kuhneao}, and Bilu-Masser-Zannier, see \cite{bilumasserzannier}. In both of the aforementioned papers, an effective bound for the number of CM points in curves in $Y(1)^2$ is established unconditionally. From the effective results of K\"uhne and Bilu-Masser-Zannier, the assertion in \Cref{corollaryintro} would follow, without any restriction to the curve $C$, thanks to work of F. Breuer, see Theorem $4.1$ in \cite{breuer}.

\begin{remark}
We finally note that using Andr\'e's G-functions method, in fact through Andr\'e's original Theorem $1.3$ in Chapter $X$ of \cite{andre1989g}, Binyamini-Masser have announced in \cite{binyaminimasser} effective results of Andr\'e-Oort-type in $\mathcal{A}_g$. 
\end{remark}

\subsection{Summary}

We start in \Cref{section:background} with some general background on Andr\'e's G-functions method. The main result here, \Cref{gfunsthm},  encodes in a sense the interplay between G-functions and relative periods of the variation of Hodge structures given by $R^1f_{*}\Q$, where $f:\CX=\CE_1\times\ldots\times \CE_n\rightarrow S$ is some $1$-parameter family of products of elliptic curves. The main technical parts are heavily based on recent work on the G-functions method, mainly the exposition of \cite{daworr4,daworr5,papasbigboi,papaszp,davidg}. At the end of the day, given a $1$-parameter family over a number field as above, we can associate to it a naturally defined family of G-functions which we denote by $\mathcal{Y}$.

In \Cref{section:trivial} based on our previous work in \cite{papasbigboi}, mainly $\S$ $7$ there, we give a description of the so called ``trivial relations'' among the G-functions in our family. This is achieved working as in \cite{papasbigboi} via a monodromy argument using the Theorem of the Fixed Part of Andr\'e, see \cite{andrefixed}.

We continue in \Cref{section:cm}, which is the start of the main technical part of our exposition. Here we construct relations among the archimedean values of our family $\mathcal{Y}$ of G-functions at, essentially, points $s\in S(\bar{\Q})$ over which the fiber of the morphism $f$ above reflects an intersection of the image of our curve in the moduli space $Y(1)^n$ with a CM point. 

We conclude the main technical part of this text in \Cref{section:heightbounds} by establishing the height bounds needed to deduce our Large Galois orbits statements. To do this it is crucial that we assume that the abelian scheme in question ``degenerates'', namely that there exists some curve $S'$ with $S\subset S'$ and some point $s_0\in S'(\bar{\Q})$ such that the fiber at $s_0$ of the connected N\'eron model $\CX'$ of $\CX$ over $S'$ has some $\G_m$ component. The proof then is done by essentially appealing to the ``Hasse Principle'' of Andr\'e-Bombieri for the values of G-functions. To do this we show that the relations constructed in the previous sections among the values of our G-functions at points of interest are both ``non-trivial'', i.e. they do not hold generically, and ``global'', i.e. they hold for all places with respect to which our point of interest $s$ is ``close'' to the point of degeneration $s_0$. This final step, i.e. the globality of our relations, is achieved by an analogue of the original argument of Andr\'e in \cite{andre1989g} making use of Gabber's lemma to show that the points we are considering cannot be ``close'' to $s_0$ with respect to any finite place.

The reduction of \Cref{maintheorem} to the aforementioned height bounds is well-known, using so-called ``endomorphism estimates'' of Masser-W\"ustholz \cite{mwendoesti}. We include this standard argument for reasons of completeness in \Cref{section:applications}.\\

\subsection{Notation}\label{section:notation}
We introduce some notation that we adopt throughout the text.\\

Given a number field $L$ we write $\Sigma_L$ for the places of $L$, $\Sigma_{L,\infty}$ for the set of its archimedean places, and respectively $\Sigma_{L,f}$ for the set of its finite places. Then given a place $v\in \Sigma_L$ we write $\C_v$ for the complete, with respect to $v$, algebraically closed field corresponding to the place $v$. We will also write $\iota_v:L\hookrightarrow \C_v$ for the embedding of $L$ in $\C_v$ that corresponds to $v$.

Given a scheme $U$ defined over $L$, where $L$ is either a number field or $L=\bar{\Q}$, and $\iota:L\hookrightarrow \C$ an embedding of $L$ into $\C$, we write $U_{\iota}:=U\times_{L,\iota}\C$ for the base change of $U$ over $\C$.

Consider a power series $y:=\Sum{i=0}{\infty}y_ix^i\in L[[x]]$, with $L$ a number field, and let $\iota_v$ be as above the embedding that corresponds to some place $v\in \Sigma_{L}$. We write $\iota_v(y(x))$ for the power series $\Sum{i=0}{\infty}\iota_v(y_i)x^i\in\C_v[[x]]$.

Finally, for a family of such power series $y_j\in L[[x]]$ and an embedding $\iota_v:L\hookrightarrow \C_v$, we define $R_v(\{y_1\ld y_N\}):=\max\{ R_v(\iota_v(y_j))\}$, where $R_v(f)$ for a power series $f\in \C_v[[x]]$ denotes the radius of convergence of $f$. 
	

\section{Recollections on the G-functions method}\label{section:background}

The main object of study in this paper is essentially the transcendence properties of values of certain G-functions that appear either as relative periods of $1$-parameter families of products of elliptic curves or are closely related to those in a manner that we soon make specific. In this first section we review this relation in this context. 

\subsection{Our setting}\label{section:setting}

Instead of working with a curve $C\subset X(1)^n$ in the majority of our exposition we will deal with a slightly different setting modeled towards applying Andr\'e's G-function method. We dedicate this subsection to recalling this setup and the main conventions we make.\\

We consider $S'$ a smooth, not necessarily projective, geometrically irreducible curve defined over a number field $K$, a point $s_0\in S'(K)$, and set $S\subset S'\backslash \{s_0\}$ an open subset. We also assume that we are given an abelian scheme of the form $f:\CX=\CE_1\times \ldots\times\CE_n\rightarrow S$, where for each $1\leq i\leq n$ the morphism $f_i:\CE_i\rightarrow S$ defines an elliptic curve over $S$, the morphism also being defined over $K$.

Up to replacing $S'$ by a finite cover, we may assume that for each $1\leq i\leq n$ the connected N\'eron model $f_i':\CE'_i\rightarrow S'$ of the elliptic curve $f_i:\CE_i\rightarrow S$ exists. By the term ``connected N\'eron model'', we mean the group subscheme of the N\'eron model $N_{S'}(\CE_i)$ of $\CE_i$ over $S'$, whose fiber over every point of $S'$ is the connected component of the identity in the respective fiber of $N_{S'}(\CE_i)$. For more details on this we point the interested reader to $\S 2.1$ in \cite{papaszpiny1}.

In what follows, we will denote the fiber product over $S'$ of the above morphisms by \begin{center}
	$f':\CX':=\CE'_1\times\ldots \times \CE'_n\rightarrow S'$.
\end{center}Note that $\CX'$ will also be the connected N\'eron model of $\CX$ over $S'$ by standard properties of N\'eron models, see for example \cite{neron}.

With \Cref{defcoord} in mind we introduce the following:
\begin{definition}\label{defsmorsingcoordrel}Let $S'$, $s_0$, and $f'$ be as above. The coordinate $i$ is said to be smooth for the pair $(S',s_0)$ if $(f'_i)^{-1}(s_0)$ is an elliptic curve. A smooth coordinate $i$ for the pair $(S',s_0)$ will be called a CM coordinate for $(S',s_0)$ if in addition $(f'_i)^{-1}(s_0)$ is a CM elliptic curve. On the other hand, the coordinate $i$ is said to be singular for our pair $(S',s_0)$ if it is not smooth, i.e. if $(f'_i)^{-1}(s_0)\simeq \G_m$.\end{definition}

Upon passing to a finite \'etale cover of $S'$ we may assume, which we do for the rest of this section without further mention, that the following assumption holds for our setting:
\begin{assumption}\label{assummonodromy}
The local monodromy around $s_0$ acts unipotently on the fibers of $R^1(f_k)_{*}\Q$ in some analytic neighborhood of $s_0$, for all singular coordinates $k$ for $(S',s_0)$.
\end{assumption}
\begin{remark}
	This is essentially a usage of Schmid's results on the quasi-unipotent action of local monodromy, see \cite{schmid}. Such an assumption does not interfere with our intended purpose of bounding the height of points $s\in S'(\bar{\Q})$ for which the fiber of $f$ is a tuple of CM elliptic curves. For more on this we point the interested reader to the proof of Lemma $2.4$ in \cite{papaszpiny1}.
\end{remark}
\subsubsection{Relative periods}\label{section:periodsbasics}

As was already mentioned, our main objects of study will be relative periods of morphisms and their arithmetic properties. Here we have tried to summarize the construction of these objects in the particular case we study, i.e. that of a product of elliptic curves over a curve. For a more general account we point the interested reader to Chapter $IX$ of \cite{andre1989g}, whose exposition and notation has clearly influenced ours.

We follow our earlier notation from this section and furthermore fix a place $v\in \Sigma_{K,\infty}$ with corresponding embedding $\iota_v:K\hookrightarrow \C$. We then have a canonical isomorphism, due to A. Grothendieck, between the de Rham and singular cohomologies on the relative level \begin{equation}\label{eq:grothver1}
	H^1_{DR}(\CX/S)\otimes_{\CO_{S}}\CO_{S^{an}_v}\rightarrow R^1(f_v)_{*}(\Q)\otimes_{\Q}\CO_{S^{an}_v},
\end{equation}where $S^{an}_v$ stands for the analytification of $S_v$.

In our particular situation, i.e. that  of an $n$-tuple of elliptic curves over $S$, we will furthermore have that\begin{equation}\label{eq:grothver2}
		H^1_{DR}(\CX/S)=H^1_{DR}(\CE_1/S)\oplus\ldots\oplus H^1_{DR}(\CE_n/S),\text{ and}
	\end{equation}
	\begin{equation}\label{eq:grothsplitbetti}
		R^1(f_v)_{*}(\Q)=R^1(f_{1,v})_{*}(\Q)(1)\oplus \ldots\oplus R^1(f_{n,v})_{*}(\Q)(1).
\end{equation}where we think of $R^1(f_{k,v})_{*}\Q(1)$ as the variation of Hodge structures whose fibers are the first cohomology group of the corresponding fibers of $f_k$. We also note that the isomorphism \eqref{eq:grothver1} is compatible with the direct sum decompositions in \eqref{eq:grothver2} and \eqref{eq:grothsplitbetti}.

Let us choose for each $1\leq k\leq n$ a basis of sections $\{\omega_{2k-1},\omega_{2k}\}$ of $H^1_{DR}(\CE_k/S)|_{U}$ over some affine open $U$, a trivializing frame $\Gamma_{k,v}=\{\gamma_{2k-1,v},\gamma_{2k,v}\}$ of $R^1(f_{k,v})_{*}(\Q)(1)|_{V}$ over some simply connected $V\subset U_{v}$, and set $\Gamma_v:=\Gamma_{1,v}\sqcup\ldots\sqcup\Gamma_{n,v}$ which will be a trivializing frame of the local system $(R^1(f_{1,v})_{*}(\Q)(1)\oplus \ldots\oplus R^1(f_{n,v})_{*}(\Q)(1))^{\vee}|_V$. 

For each $k$, associated to the above, we then get a matrix of relative periods of $V$ which we denote by \begin{equation}\label{eq:relpercoord}
\CP_{\Gamma_{k,v}}:=\begin{pmatrix}
\frac{1}{2\pi i}\int_{\gamma_{2k-1,v}} \omega_{2k-1}   &  \frac{1}{2\pi i}\int_{\gamma_{2k,v}}\omega_{2k-1}   \\
\frac{1}{2\pi i}\int_{\gamma_{2k-1,v}}  \omega_{2k}  & \frac{1}{2\pi i}\int_{\gamma_{2k,v}}\omega_{2k}
\end{pmatrix}
\end{equation}which encodes the restriction of the canonical isomorphism $H^1_{DR}(\CE_k/S)\otimes\CO_{S^{an}_v}\rightarrow (R^1(f_{k,v})_{*}\Q(1))^{\vee}\otimes\CO_{S^{an}_v}$ over the open analytic set $V$.

Similarly, associated to the chosen basis $\{\omega_i:1\leq i\leq 2n\}$ and the trivializing frame $\Gamma_v$ as above, we get a matrix of relative periods encoding the isomorphism \eqref{eq:grothver2} which we will denote by $\CP_{\Gamma_v}$. We note that by construction of our trivializing frame and basis $\{\omega_i\}$ this matrix will be block diagonal, since the isomorphism in question respects the splitting in de Rham and Betti cohomology given by $\CX=\CE_1\times\ldots\times\CE_n$, and the diagonal blocks will be the matrices $\CP_{\Gamma_{k,v}}$ above. 

\begin{remark}
	We have opted for a notation that does not mention either the choice of a basis or that of a simply connected $V$ over which we get a trivializing frame. The reason for that is that pretty much throughout this text we will consider a fixed such basis $\omega_i$, appropriately chosen, and care more to encode the family of relative periods that come out of \eqref{eq:grothver2} as one varies the chosen place $v\in \Sigma_{K,\infty}$. 
\end{remark}

	\subsection{G-functions and relative periods}\label{section:gfuns}

In this subsection we momentarily abandon the setting in \Cref{section:setting} that we adopt almost throughout the rest of this text. Namely for now, we consider a fixed $f':\CX'\rightarrow S'$, this time defined over $\bar{\Q}$, $s_0\in S'(\bar{\Q})$ with the same properties as in \Cref{section:setting}, and an embedding $\iota:\bar{\Q}\rightarrow \C$. Throughout this subsection we also fix a local parameter $x$ of $S'$ at $s_0$. Later on, see \Cref{section:goodcovers}, we will be more careful about this choice when we review what we call a ``good cover of the curve $S'$''.

The main objects of our study will be G-functions that are naturally associated to a pair $(\CX'\rightarrow S',s_0)$ as above. These objects were first studied by C. L. Siegel in \cite{siegel}. For the convenience of the reader we include here the following:
\begin{definition}
	Let $K$ be a number field and let $y=\Sum{n=0}{\infty}a_n x^n\in K[[x]]$. Then $y$ is called a \textbf{G-function} if the following are true: 
	\begin{enumerate} 	
		\item $\forall v\in \Sigma_{K,\infty}$ we have that $i_v(y)\in \C_v[[x]]$ defines an analytic function around $0$,
		\item there exists a sequence $(d_n)_{n\in\N}$ of natural numbers such that 
		\begin{itemize}
			\item $d_n a_m\in \mathcal{O}_K$ for all $m\leq n$,
			\item there exists $C>0$ such that $d_n\leq C^n$ for all $n\in \N$,
		\end{itemize}	
		\item $y$ satisfies a linear homogeneous differential equation with coefficients in $K(x)$.
	\end{enumerate}
\end{definition}

We also introduce the following convenient:
\begin{definition}\label{defgmatrix}
	We call a matrix $A\in M_{r_1\times r_2}(\bar{\Q}[[x]])$ a \textbf{G-matrix} if all of its entries are G-functions.
\end{definition}

\begin{theorem}\label{gfunsthm}
	There exists a basis of sections $\{\omega_i:1\leq i\leq 2n\}$ of $H^1_{DR}(\CX/S)$ over $U:=U'\backslash \{s_0\}$, where $U'$ is some open affine neighborhood of $s_0$, and an associated family of G-matrices $Y_{G,k}=(y_{i,j,k})\in \GL_2(\bar{\Q}[[x]])$ such that, writing $\mathcal{Y}:=\{y_{i,j,k}:1\leq i,j\leq 2, 1\leq k\leq n\}$, for every $s\in U(\bar{\Q})$ with $|x(s)|_\iota<\min\{1,R_\iota(\mathcal{Y})\}$ we have that\begin{enumerate}
		\item if $k$ is a smooth coordinate for $(S',s_0)$ then there exists a symplectic trivializing frame $\Gamma_{k,\iota}=\{\gamma_{2k-1,\iota},\gamma_{2k,\iota}\}$ of $(R^1(f_{k,\iota})_{*}(\C))^{\vee}|_V$ over some small enough analytic neighborhood $V\subset S_\iota$ of $s$ such that \begin{equation}\label{eq:gfunssmooth}
		\CP_{\Gamma_{k,\iota}}(s)=\iota(Y_{G,k}(x(s)))\cdot \Pi_{k,\iota},
		\end{equation}where $\Pi_{k,\iota}\in \GL_2(\C)$ is such that, if the coordinate $k$ is furthermore CM for $(S',s_0)$, it is of the form \begin{equation}
		\Pi_{k,\iota}=\begin{pmatrix}
			\frac{\varpi_{k,\iota}}{2\pi i}&0\\0&\varpi_{k,\iota}^{-1}.
		\end{pmatrix}
	\end{equation}

		\item if $k$ is a singular coordinate for $(S',s_0)$ there exist $d_k$, $d'_k\in \bar{\Q}$ independent of the chosen embedding $\iota$, and a symplectic trivializing frame $\Gamma_{k,\iota}=\{\gamma_{2k-1,\iota},\gamma_{2k,\iota}\}$ of $(R^1(f_{k,\iota})_{*}\Q)^{\vee}|_V$ over some small enough analytic neighborhood $V\subset S_\iota$ of $s$ such that \begin{equation}\label{eq:gfunssmooth}
		\CP_{\Gamma_{k,\iota}}=	\iota(Y_{G,k}(x(s)))\cdot \Pi_{k,\iota} \cdot \begin{pmatrix}
			1&N_k\log \iota(x(s))\\0&1
		\end{pmatrix},
		\end{equation}where $N_k\in \Q$ and $\Pi_{k,\iota}\in\GL_2(\C)$ is such that its first column is $\begin{pmatrix}\iota(d_k)\\ \iota(d'_k)\end{pmatrix}$.
	\end{enumerate}
\end{theorem}

\begin{remarks}$1$. We stress that the choices of the bases and the various trivializations in the previous theorem are independent of the point $s\in S(\bar{\Q})$ in question but depend on the ``base'' point $s_0$. The various frames will also obviously depend on the choice of the chosen embedding $\iota$. We return to this last dependence in the next subsection.\\
	
	$2$. From the previous theorem and the remarks in \Cref{section:periodsbasics} we know that the relative period matrix $\CP_{\Gamma_\iota}$ associated to the morphism $f:\CX\rightarrow S$, the embedding $\iota$, the basis $\{\omega_i:1\leq i\leq 2n\}$, and the frame $\Gamma_\iota=\Gamma_{1,\iota}\sqcup\ldots\sqcup\Gamma_{n,\iota}$ will be block diagonal with diagonal blocks the above $\CP_{\Gamma_{k,\iota}}$ which are described as in \Cref{gfunsthm}.\\
	
	$3$. We expect that this result is known to experts in the area. Indeed the ideas here appear already in \cite{andre1989g} and \cite{andremots} though the theorem itself is not expressly stated in this format.	
\end{remarks}

We start with the following fundamental lemma about periods of CM elliptic curves that we will need in the proof of the above theorem.
\begin{lemma}\label{lemmacmperiods}
	Let $E/L$ be a CM elliptic curve defined over a number field $L$ and assume that $F:=\End^{0}_{\bar{\Q}}(E)=\End^{0}_L(E)$. We fix an embedding $\iota_v:L\hookrightarrow\C$, corresponding to some $v\in \Sigma_{L,\infty}$ and let $V_{dR}:=H^1_{DR}(E/L)$ and $V_{\Q}:=H_1(E_{v},\Q)$.

	Then there exist \begin{enumerate}
		\item a symplectic basis $\omega_1$, $\omega_2$ of $V_{dR}\otimes_L LF$, and 
		\item a symplectic basis $\gamma_1$, $\gamma_2$ of $V_{\Q}\otimes LF$, 
	\end{enumerate}
such that the period matrix of $E$ with respect to these choices is of the form \begin{equation}\begin{pmatrix}
		\frac{\varpi_v}{2\pi i}&0\\0&\varpi_v^{-1},
	\end{pmatrix}
\end{equation}for some $\varpi_v\in \C$.\end{lemma}
	\begin{proof}Via the action of $F$ on $V_{dR}$ and $V_{\Q}$ we get splittings of $V_{dR,LF}$ and $V_{\Q,LF}$ which are compatible via Grothendieck's comparison isomorphism\begin{center}
		$P:V_{dR}\otimes_{L}\C\rightarrow (V_{\Q})^{\vee}\otimes_{\Q} \C$.
	\end{center}
	
	In more detail, on the one hand we have the splitting on the de Rham side:\begin{equation}\label{eq:splitderham}
		V_{dR}\otimes_{L} LF =W^{\sigma_1}_{dR}\oplus W^{\sigma_2}_{dR}, 
	\end{equation}and the splitting on the Betti side:\begin{equation}\label{eq:splitbetti}
		V_{\Q}\otimes LF= W_{\sigma_1}\oplus W_{\sigma_2},
	\end{equation}where $\sigma_i:F\hookrightarrow \C$ are the two embeddings of $F$ in $\C$. Here, following the notation in \cite{andre1989g} Ch. X, we denote by $W_\sigma$ and $W^{\sigma}_{dR}$ the subspaces of the respective vector space where $F$ acts via the embedding $\sigma:F\rightarrow \C$.
	
	By Lemma $8.2 $ of \cite{papasbigboi}, also its ``dual'', we have that there exist the following:\begin{enumerate}
		\item a symplectic basis $\omega_1$, $\omega_2$ of $V_{dR ,LF}$  for which we furthermore have that $\omega_i $ spans $W^{\sigma_i}_{dR}$, and 
		
		\item $\gamma_1$, $\gamma_2$ a symplectic basis of $V_{\Q,LF}$ such that $\gamma_i$ spans the subspace $W_{\sigma_i}$.
	\end{enumerate}
	
Note that we have \begin{equation}\label{eq:cmeq1}
		P(\omega_i)=\left( \frac{1}{2\pi i}\int_{\gamma_1}^{}\omega_i\right) \gamma_1^{\vee}+\left(\frac{1}{2\pi i}\int_{\gamma_2}^{}\omega_i\right) \gamma_2^{\vee}, \text{ } i=1,2.
	\end{equation}One then has from the compatibility of the action of $F$ with this isomorphism, that for every $\lambda\in F$:\begin{equation}\label{eq:cmeq2}
		P(\lambda \omega_i)=\left(\frac{1}{2\pi i}\int_{\gamma_1}^{}\omega_i\right) \sigma_1(\lambda)\gamma_1^{\vee}+\left(\frac{1}{2\pi i}\int_{\gamma_2}^{}\omega_i\right) \sigma_2(\lambda)\gamma_2^{\vee}, \text{ } i=1,2.
	\end{equation}
	
	On the other hand we have from the definition of the $\omega_i$ that \begin{equation}\label{eq:cmeq3}P(\lambda\omega_i)=\sigma_i(\lambda)P(\omega_i).
	\end{equation}Since all of the above is true for any $\lambda\in F$, by comparing coefficients with \eqref{eq:cmeq2} we get \begin{equation}\label{eq:cmeq4}
		\frac{1}{2\pi i}\int_{\gamma_1}\omega_2=\frac{1}{2\pi i}\int_{\gamma_2}\omega_1=0.
	\end{equation}
	
	Now set $\varpi_1=\frac{1}{2\pi i} \int_{\gamma_1}\omega_1$ and $\varpi_2:=\frac{1}{2\pi i}\int_{\gamma_2}\omega_2$. Then the Legendre relations give $\varpi_1\cdot \varpi_2=\frac{1}{2\pi i}$. In particular we get that the period matrix with respect to these choices of bases is of the form \begin{equation}\label{eq:cmeq5}
		\begin{pmatrix}
			\frac{\varpi_v}{2\pi i}&0\\ 0& \varpi_v^{-1},
		\end{pmatrix}
	\end{equation}as we wanted.\end{proof}

\begin{remark}
	We note that for the $\varpi_v\in\C$ it is known that $tr.d._{\bar{\Q}}(\varpi_v,\pi)=2$. This follows from Grothendieck's period conjecture which is known here by work of Chudnovsky, see \cite{chudicm,chudbook}.
\end{remark}
	\begin{proof}[Proof of \Cref{gfunsthm}] Part $(2)$ is \cite{andre1989g}, Ch. IX, $\S 4$, Theorem $1$ when $g=1$. We note that the explicit description of the period matrix is inherent in the proof. See also the proof of Claim $3.7$ in \cite{papaszp} where this explicit description appears. We note that \Cref{assummonodromy} is needed here, see the proof of Theorem $3.1$ of \cite{papaszp} for more details.
	
	The matrix $Y_{G,k}(x)$ will be the normalized uniform solution of the G-operator $\vartheta-G_k$, where $\vartheta:=x\frac{d}{dx}$ and $G_k=(g_{i,j,k})$ is given by $\nabla_{\vartheta}(\omega_{i,k})=\Sum{j=1}{2} g_{i,j,k}\omega_{i,j.k}$, where $\nabla$ denotes the Gauss-Manin connection in question. The fact that $G_k$ is a G-operator\footnote{See \cite{andre1989g} page $76$ for a definition of this notion.} follows from the proof of the theorem in the appendix of Chapter V in \cite{andre1989g}, since in this case the operator corresponds to a geometric differential equation. That the entries of the matrix $Y_{G,k}$ are G-functions now follows from the Corollary in Ch. V, $\S 6.6$ of loc. cit..
	
	So for the singular coordinates we choose the basis $\{\omega_{2k-1},\omega_{2k}\}$ and a symplectic trivializing frame $\Gamma_{k,\iota}$ of $H^1_{DR}(\CE_k/S)|_U$ and $R^1(f_{k,\iota})_{*}(\Q)(1)|_V$ respectively as specified in Theorem $3.1$ of \cite{papaszp}.\\
	
	Now we move on to the proof of $(1)$ and the smooth coordinates for $(S',s_0)$. For the non-CM smooth coordinates our work is simpler. Namely we may choose any symplectic basis of $\{\omega_{2k-1},\omega_{2k}\}$ of $H^1_{DR}(\CE'_k/S')$ over some neighborhood $U'$ of $s_0$ and any symplectic frame of $R^1(f'_{k,\iota})_{*}(\C)(1)|_V$ for some small enough analytic neighborhood $V$ of $s_0$. 
	
	To see this, first of all note that in this case the differential system $\vartheta-G_k$ that arises as above is such that $G_k(0)=0$. Indeed, in this case the morphisms $f'_k:\CE'_k\rightarrow S'$ are in fact smooth and proper. Therefore, $G_k(0)$ which coincides with the residue of the connection at the point $s_0$ will be $0$.
	
	Now any solution of the system $\vartheta-G_k$ will be of the form $X_k=Y_{G,k}\cdot \Pi_{k,\iota}$ where $\Pi_k\in \GL_2(\C)$, see \cite{andre1989g} Ch. III, $\S 1$. Since $\CP_{\Gamma_{k,\iota}}$ is such a solution for any choice of $\Gamma_{k,\iota}$ we are done. We note that by construction we will also have \begin{equation}\label{eq:peridosatbase}
		\CP_{\Gamma_{k,\iota}}(0)=\Pi_{k,\iota}
	\end{equation}where $\Pi_{k,\iota}=\begin{pmatrix}
	\frac{1}{2\pi i}\int_{\gamma_{2k-1,\iota}}^{}(\omega_{2k-1})_{s_0}&\frac{1}{2\pi i}\int_{\gamma_{2k,\iota}}^{}(\omega_{2k-1})_{s_0}\\
	\frac{1}{2\pi i}\int_{\gamma_{2k-1,\iota}}^{}(\omega_{2k})_{s_0}&\frac{1}{2\pi i}\int_{\gamma_{2k,\iota}}^{}(\omega_{2k})_{s_0}
\end{pmatrix}$ will be the period matrix of the elliptic curve $(\CE_k')_{s_0}$.

Let us finally look at the CM coordinates. Using \Cref{lemmacmperiods} we can then find a symplectic basis $\{\omega_{2k-1},\omega_{2k}\}$ of $H^1_{DR}(\CE'_k/S')|_{U'}$ and a symplectic trivializing frame of the local system $R^1(f'_{k,\iota})_{*}(\C)(1)|_{V'}$ in a small enough neighborhood $V'$ of $s_0$ as above with the properties we wanted. Whence the description of $\Pi_{k,\iota}$ when $k$ is a CM coordinate follows. 

Finally, in both cases, i.e. CM or non-CM smooth coordinate, the fact that the matrix $Y_{G,k}$ is a G-matrix follows from the same exact argument as in the singular case above.	\end{proof}
	\subsubsection{Family of G-functions associated to $s_0$}

Let $f:\CX'\rightarrow S'$ be defined over $\bar{\Q}$ and $s_0\in S'(\bar{\Q})$ be as above.

Our first order of business is to associate from now on a family of G-functions to the point $s_0$. The ``natural expectation'' to associate to $s_0$ the entire family $\mathcal{Y}$ as defined in \Cref{gfunsthm} turns out to give various complications down the line. First of all, only the first column of relative periods $\CP_{\Gamma_{k,\iota}}$ with $k$ singular will play an actual role in what we need. Secondly, the so called ``trivial relations'' of the family $\mathcal{Y}$ are messier to describe. 

With these goals in mind, let us fix for now a singular coordinate $k$. Then from \Cref{gfunsthm} we know that locally near $s_0$ \begin{equation}\label{eq:firstcolumn1}
	\CP_{k,v}=\iota_v(Y_{G,k})\cdot\Pi_{k,\iota}\cdot \begin{pmatrix} 1&N_k\log(\iota(x))\\0&1\end{pmatrix}.
\end{equation}In particular for our choice, in the proof of \Cref{gfunsthm}, of basis $\omega_{2k-1}$, $\omega_{2k}$ of $H^1_{DR}(\CE_k/S)|_U$ and trivialization $\Gamma_{k,\iota}$ of the local system $R^1(f_{k,\iota})_{*}\Q(1)$ the first column of the matrix $\CP_{k,\iota}$ will be of the form \begin{equation}\label{eq:firstcolumn2}
\begin{pmatrix}
	\frac{1}{2\pi i}\int_{\gamma_{2k-1,\iota}}\omega_{2k-1}\\	\frac{1}{2\pi i}\int_{\gamma_{2k-1,\iota}}\omega_{2k}
\end{pmatrix}= \begin{pmatrix}
\iota(d_ky_{1,1,k}(x)+d'_ky_{1,2,k}(x))\\\iota(d_ky_{2,1,k}(x)+d'_ky_{2,2,k}(x)).
\end{pmatrix}
\end{equation}

\begin{lemma}\label{lemmasymplectic}
	Let $f_k:\CE_k\rightarrow S$ be a singular coordinate for some $f:\CX\rightarrow S$ as above. Then there exists a basis $\omega'_{2k-1}$, $\omega'_{2k}$ of $H^1_{DR}(\CE_k/S)|_U$, where $U=U'\backslash\{s_0\}$ for some possibly smaller affine neighborhood $U'$ of $s_0$ as before, such that \begin{enumerate}
		\item with respect to the trivializing frame $\Gamma_{k,\iota}$ chosen in \Cref{gfunsthm}, the entries of the first column of the relative period matrix $\CP_{k,\iota}$ are G-functions, and 
		
		\item the matrix of the polarization on $H^1_{DR}(\CE_k/S)|_{U}$ in terms of this basis is of the form \begin{equation}
			e_k\cdot \begin{pmatrix}
				0&1\\-1&0
			\end{pmatrix},
		\end{equation}with $e_k\in \mathcal{O}_{S'}(S)^{\times}$.
	\end{enumerate}
\end{lemma}
\begin{proof}We note that the basis $\omega_i$ chosen in the proof of \Cref{gfunsthm} is in fact the restriction on $U:=U'\backslash \{s_0\}$ of a basis, which we denote by the same notation, of the vector bundle $\mathcal{E}|_{U'}$, where $\mathcal{E}:=H^1_{DR}(\CE_k/S)^{can}$ is the canonical extension of the vector bundle $H^1_{DR}(\CE_k/S)$ to $S'$. 
	
	By the proof of Lemma $6.7$ of \cite{daworr5} there exist sections $\omega_1$, $\eta_1$ of $\mathcal{E}|_{U'}$, upon possibly replacing the original $U'$ by a smaller affine open neighborhood of $s_0$ in $S'$ and letting $U=U'\backslash \{s_0\}$ as before, such that $(\omega_1)|_U$, $(\eta_1)_{U}$ is a basis of $H^1_{DR}(\CE_k/S)|_{U}$, and part $(2)$ of the lemma holds. 
		
	Now note that we have, by construction, that there exists a matrix 
	$\begin{pmatrix}a&b\\c&d\end{pmatrix}\in M_{2\times 2}(\mathcal{O}(U'))$ such that \begin{equation}
		\omega_1=a\omega_{2k-1}+b\omega_{2k}\text{, and }\eta_1=c\omega_{2k-1}+d\omega_{2k}
	\end{equation}

With respect to the basis $\{\omega_1,\eta_1\}$ and the frame $\Gamma_{k,\iota}$ the first column of the relative period matrix is of the form \begin{equation}\label{eq:firstcolumnfin}
	\begin{pmatrix}
		\iota(a)\iota(F_1)(x)+\iota(b)\iota(F_2)(x)\\\iota(c)\iota(F_1)(x)+\iota(d)\iota(F_2)(x)
	\end{pmatrix},
\end{equation}where $F_i(x)$ are the entries of \eqref{eq:firstcolumn2}, which will be G-functions by the preceding discussion.

The lemma on page $26$ of \cite{andre1989g} and the proposition on page $27$ of loc. cit. show that the $a$, $b$, $c$, $d$ have power series expansions on $x$ that are G-functions. From Theorem $D$ in the introduction of loc. cit. the coordinates of the vector in \eqref{eq:firstcolumnfin} will be G-functions. We thus set $\omega'_{2k-1}:=\omega_1$ and $\omega_{2k}:=\eta_1$.
\end{proof}

\begin{definition}\label{defassoctoapoint1}
	We denote by $\mathcal{Y}_{s_0}$ the family of G-functions that consists of the following power series:\begin{enumerate}
		\item the entries of the G-matrices $Y_{G,k}:=(y_{i,j,k}(x))$ appearing in \Cref{gfunsthm} for all smooth coordinates $k$ of $S'$, and 
		
		\item the entries of the first column, which we denote by $\begin{pmatrix}
			y_{1,1,k}(x)\\y_{2,1,k}(x)
		\end{pmatrix}$, of the relative period matrices $\CP_{\Gamma_{k,\iota}}$ with respect to the bases of \Cref{lemmasymplectic}.
	\end{enumerate}

We call this the\textbf{ family of G-functions associated locally to the point} $s_0$.
\end{definition}
	\subsubsection{Independence from archimedean embedding}\label{section:independent}

Let us return to our original notation with $f':\CX'\rightarrow S'$ defined over some number field $K$, $s_0\in S'(K)$, as in \Cref{section:setting} satisfying \Cref{assummonodromy}. Let us also fix for now a local parameter $x$ of $S'$ at $s_0$. 

Let $\{\omega_i:1\leq i\leq 2n\}$ be the basis of $H^1_{DR}(\CX/S)$ appearing in \Cref{gfunsthm} with the $\omega_i$ that correspond to singular coordinates replaced by the $\omega_i'$ of \Cref{lemmasymplectic}. From \Cref{gfunsthm} and \Cref{lemmasymplectic}, we then know that, upon fixing an embedding $\iota:\bar{\Q}\hookrightarrow\C$, G-functions appear in a specific way in the description of the relative periods of $f_{\bar{\Q}}:\CX_{\bar{\Q}}\rightarrow S_{\bar{\Q}}$ close to the point $s_0$. Since these G-functions are solutions to various geometric differential equations the field generated by their coefficients over $\Q$ is in fact a number field. Let us denote this field by $K_{\mathcal{Y}}$. 

We define the number field $K_{f'}$ to be the compositum of the following fields:\begin{enumerate}
	
	\item the field $K$ over which our setup is defined, 
	
	\item the CM-fields $F_k$ associated to the CM coordinates of $S'$,
	
	\item the number fields $\Q(d_k,d'_k)$ generated by the constants $d_k$, $d_k'\in \bar{\Q}$ associated themselves to each singular coordinate of the curve $S'$, and 
	
	\item the number field $K_{\mathcal{Y}}$.
\end{enumerate}

Upon base changing the morphisms $f':\CX'\rightarrow S'$ by $K_{f'}$, in essence replacing $K$ by $K_{f'}$, we may work, which we do from now on, under the following assumption:\begin{assumption}\label{assumptionfield}
	In the above setting we have $K_{f'}=K$ so that all the constants that appear in \Cref{gfunsthm} associated to the relative periods of $f$ near $s_0$ are in fact in the base number field $K$.
\end{assumption}

For every archimedean embedding $\iota_v:K\hookrightarrow \C$, associated to an archimedean place $v\in\Sigma_{K,\infty}$, we may repeat the process of \Cref{gfunsthm} and \Cref{lemmasymplectic}, keeping the basis $\omega_i$ of $H^1_{DR}(\CX/S)|_{U}$ chosen for a fixed place $v_0\in \Sigma_{K,\infty}$. It is easy to see that all the algebraic constants, i.e. the coefficients of the G-matrices and the $d_k$, $d'_k$, depend only on the choice of that basis. One can then find trivializing frames of $R^1(f_{k,v})_{*}(\C)$ for the various coordinates $k$, with $v\in \Sigma_{K,\infty}$ and $v\neq v_0$, such that the relative periods of the morphism $f$ are of the form described in \Cref{gfunsthm}. The only non-trivial case, that of singular coordinates, is dealt with by the lemma in Ch. X, $\S 3.1$ of \cite{andre1989g}. 

In other words we have the following
\begin{lemma}\label{independence}
	Let $s\in S(L)$ with $L/K$ finite and let $v\in \Sigma_{L,\infty}$ be such that $|x(s)|_v<\min \{1, R_v(\mathcal{Y}_{s_0})\}$. Then there exists a choice of a trivializing frame $\Gamma_{v}$ of $(R^1(f_{1,v})_{*}(\Q)(1)\oplus \ldots\oplus R^1(f_{n,v})_{*}(\Q)(1))^{\vee}|_V$ for some small enough analytic neighborhood $V$ of $s$ in $S_v$ such that \begin{enumerate}
\item $\CP_{k,v}(s)=\iota_v(Y_{G,k}(x(s)))\cdot \Pi_{k,v}$ for all smooth coordinates $k$ of $S$, and 

\item the first two columns of the relative period matrix $\CP_{k,v}(s)$ are $\begin{pmatrix}
	\iota_v(y_{1,1,k}(x(s)))\\	\iota_v(y_{2,1,k}(x(s)))
\end{pmatrix}$, for all singular coordinates $k$ of $S'$.
	\end{enumerate}
\end{lemma}
	\subsection{Height bounds}\label{section:goodcovers}

We consider a fixed projective curve $S'$ and associated semiabelian scheme $f':\CX'=\CE'_1\times\ldots\times \CE'_n\rightarrow S'$ defined over a number field $K$ as in \Cref{section:setting}. As usual we also fix a point $s_0\in S'(K)$ which is such that $\CX'_{s_0}$ is not an abelian variety and $S\subset S'\backslash \{s_0\}$ such that $f=f'|_S$ is a product of families of elliptic curves. We assume from now on that \Cref{assummonodromy} and \Cref{assumptionfield} hold for $\CX'\rightarrow S'$ as well as the following:
\begin{assumption}\label{hodgegenericity}
	The image $m(S)$ of $S$ via the morphism $m:S\rightarrow Y(1)^n$, which is induced from the scheme $f:\CX\rightarrow S$, is a Hodge generic curve. In other words, $m(S)$ is not contained in a proper special subvariety of $Y(1)^n$, as in the setting of \Cref{maintheorem} and \Cref{corollaryintro}.
\end{assumption}

\begin{definition}\label{defgadmi}We say that the semiabelian scheme $\CX'\rightarrow S'$ is $G_{AO}$-admissible with respect to $s_0$ if all of the above hold and furthermore either of the following holds:\begin{enumerate}
		\item there exists at least one CM coordinate for $(S',s_0)$, or
		
		\item there exist at least two singular coordinates for $(S',s_0)$. \end{enumerate}\end{definition}

The output we need from the G-functions method in order to establish \Cref{maintheorem} is the following:

\begin{theorem}\label{heightboundandreoort} 
	Let $f':\CX'\rightarrow S'$ be a $G_{AO}$-admissible semiabelian scheme with respect to some $s_0\in S'(K)$, for which $\CX'_{s_0}$ is not an abelian variety. Then there exist effectively computable constants $c_1$ and $c_2$, depending on $f'$, $S'$, and $s_0$, such that for all $s\in S'(\bar{\Q})$ for which the fiber $\CX_s$ is CM we have that \begin{equation}\label{eq:htbdao}
		h(s)\leq c_1[K(s):\Q]^{c_2}.
	\end{equation}	
\end{theorem}

Following the exposition of \cite{daworr4,daworr5} we reduce here the proof of \Cref{heightboundandreoort} to the following analogue of Lemma $6.6$ of \cite{daworr5}:

\begin{lemma}\label{htboundred}It suffices to prove \Cref{heightboundandreoort} under the additional assumptions that there exists a regular $\CO_K$-model $\mathfrak{S}$ of $S'$, a semiabelian scheme $\mathfrak{X}\rightarrow \mathfrak{S}$, and a rational function $x\in K(S')$ such that: 
	\begin{enumerate}[label=(\roman*)]
		\item $\mathfrak{X}_K \simeq \CX'$ as semiabelian schemes over $S'$,
		
		\item all zeroes of $x$ are simple and lie in $S'(K)$,
		
		\item the induced cover $x:S'\rightarrow \mathbb{P}^1$ is Galois, i.e. the group $\aut_{x}(S') $ of automorphisms of $S'$ that commute with $x$ acts transitively on the $\bar{\Q}$-fibers of $x$, 
		
		\item\label{condition3} for each $1\leq i\leq n$ the $i$-th coordinate is of the same type (i.e. singular, smooth, CM) for all pairs $(S',\xi)$ where $\xi$ is a root of $x$, and
		
		\item\label{condition4} the abelian part of the fiber $\CX'_{\xi}$ has everywhere semistable reduction for all roots $\xi$ of $x$.
	\end{enumerate}
	\end{lemma}
	\begin{proof}The idea is to find a finite cover of our original curve $S'$ that satisfies all of the above and to establish \Cref{heightboundandreoort} there. In order to do that we may have to replace $K$ by a finite extension, which is harmless as far as the validity of \Cref{heightboundandreoort} is concerned.
		
The proof of the existence of the cover of $S'$ with the required properties is identical to that of Lemma $6.6$ in \cite{daworr5}. The only thing we need to establish here is that the semiabelian scheme $\mathfrak{X}$ in question in fact splits as a product of semiabelian $1$-dimensional schemes so that \ref{condition3} and \ref{condition4} above hold for its generic fiber over $\spec(\CO_K)$. 
		
For \ref{condition4} we note that by replacing the original morphism $f':\CX'\rightarrow S'$, that is $G\rightarrow C$ in the notation of \cite{daworr5}, by its base change by a finite extension of $K$ we may assume due to Grothendieck's semistable reduction theorem, see Th\'eor\`eme $3.6$ in \cite{sga71}, that the fiber over $s_0$ has everywhere semistable reduction. We note that in this case we know that any potential CM factors that appear in $\CX'_{s_0}$ will have everywhere good reduction. 

Following the proof of \cite{daworr5} and its notation, we know that the semiabelian scheme $\mathfrak{X}:=\mathfrak{G}_4$ will be the pullback to $\mathfrak{C}_4$, in the notation of loc. cit., so that over $C_4^{*}$, the analogue of our ``$S$'', the splitting of $\mathfrak{X}_K$ into elliptic curves holds by virtue of the general construction in loc. cit.. From this we get the required splitting over all of $C_4$ as a product of $1$-dimensional semiabelian schemes.

Finally, \ref{condition3} holds again by the construction in \cite{daworr5}, since, in their notation, the roots $\xi$ of $x$ will be such that $\nu(\xi)=\tilde{s}_0$ and hence the fiber $(\mathfrak{X}_K)_\xi$ can be written as a product of $\mathbb{G}_m$-terms and elliptic curves appearing in the exact same order as the terms of $\CX'_{s_0}$.	\end{proof}

\begin{remark}\label{remarkoneffectivityofreduction}
With the effectivity of the constants in \Cref{heightboundandreoort} in mind, it is important to note that the construction of the cover, as well as its main outputs, in \Cref{htboundred} is effective in terms of our original data, i.e. the morphism $f':\CX'\rightarrow S'$. This follows, for the most part, by the proof of Lemma $5.1$ in \cite{daworr4} and Lemma $6.6$ of \cite{daworr5}. 

The most crucial quantity we will need a good bound for is the number $l$ of zeroes of the rational function $x$. This is nothing but the number ``$d$'' that appears in the proof of Lemma $5.1$ in \cite{daworr4} and can be bounded, due to its construction via the Riemann-Roch theorem, by the genus of the curve $S'$.

The remaining crucial data we need to keep track of is the field extension of our original number field $K$, with respect to which we consider a base change of our picture. A first extension of $K$ is the one highlighted in the proof of Lemma $5.1$ in \cite{daworr4}. This will depend on our original curve $S'$. The second extension of $K$ is the one needed for condition \ref{condition4} in \Cref{htboundred}. This will depend on the fibers $\CE'_{j,s_0}$ which are CM elliptic curves. In other words, on the morphism $f'$ and the point $s_0$.
\end{remark}
    \subsubsection{Working with a good model}\label{section:gcsetting}

Let us fix therefore a semiabelian scheme $f':\CX\rightarrow S'$ defined over some number field $K$, together with rational function $x\in K(S')$ satisfying the properties outlined in \Cref{htboundred}. It is in this setting that we focus on from now on, largely following the discussion in $\S6. E$ of \cite{daworr5}. Let us denote by $\{\xi_1\ld\xi_l\}\subset S'(K)$ the zeroes of $x$. 

By the Galois property of the cover $x:S'\rightarrow \mathbb{P}^1$ we know that there exists for each $1\leq t\leq l$ an automorphism $\sigma_t\in \aut_{x}(S')$ such that $\sigma_t(\xi_1)=\xi_t$. Letting $S\subset S'\backslash\{\xi_1,\ldots,\xi_l\}$ be an affine open such that $f'|_{S}$ is smooth and $S_1:=S'\cup \{\xi_1\}$ we get for each $1\leq t\leq l$, by pulling back $\CX'$ via $S_1\hookrightarrow S'\xrightarrow{\sigma_t}S'$, a semiabelian scheme $\CX'_t$ on $S_1$. 

The idea that first appeared in \cite{daworr4} is to work with the G-functions associated to $\xi_1$ with respect to a subset of the family of semiabelian schemes $\CX'_t$ that appear here. Namely, one considers a subset $\Lambda\subset \{1,\ldots,l\}$ that is a set of representatives of the equivalence relation \begin{center}
	$t\sim t'$ if $(\CX'_t)_{\eta}$ is isogenous to $(\CX'_{t'})_{\eta}$,
\end{center}where $\eta$ denotes the generic point of $S_1$. 

From now on we write $\lambda$ for an arbitrary element in $\Lambda$ and we identify, for ease of notation, $\lambda$ with the smallest element of $\{1,\ldots,l\}$ that is in the equivalence class of $\lambda$. By this abuse of notation, we may write $\CX'_{\lambda}:=\CX'_t$ for the semiabelian scheme as above for which $t$ is the minimal one for which $t\in \lambda$. For each $\lambda\in \Lambda$, as in the discussion preceding \Cref{defassoctoapoint1}, we thus get a family of G-functions associated to our point $\xi_1$ and the semiabelian $f'_{\lambda}:\CX'_\lambda\rightarrow S_1$.

\begin{definition}\label{definassgfunsglob}
	Let $f':\CX'\rightarrow S'$ be as above. We call the collection of G-functions $\mathcal{Y}:=\{\mathcal{Y}_{\lambda}:\lambda\in \Lambda\}$ the \textbf{family of G-functions associated to }$f'$.
\end{definition}

\begin{remark}\label{gaoadmremark}
	Throughout the rest of this article we will be working with the setting summarized in this subsection. In other words, we will have a pair $(\CX'\rightarrow S',\{\xi_1\ld\xi_l\})$ comprising of a semiabelian scheme over a curve $S'$ together with a distinguished set of points in $S'$, rather than a single point $s_0$. The Galois property in \Cref{htboundred} makes it so that the fibers of $\CX'\rightarrow S'$ over the points $\xi$ of our distinguished set are in fact isomorphic.
	
	With brevity and the above comment in mind, we have opted to call ``$G_{AO}$-admissible'' for the remainder of this text schemes as above that also satisfy the conditions in \Cref{defgadmi} with respect to any point of the aforementioned distinguished set, especially since we consider the latter set a part of our input.
\end{remark}

\subsubsection{$v$-adic proximity}\label{section:proximity}

The main objective of the G-functions method is the construction of ``global'' and ``non-trivial'' relations, as defined in \cite{andre1989g} Ch. $VII$, $\S 5$, among the values of families of G-functions at points of interest. In our case, these will be the values of the family $\mathcal{Y}$ of \Cref{definassgfunsglob} at points over which the fibers of $f'$ are CM. In order to do this we will need to speak of proximity of our points of interest to the roots $\xi$.

Given a point $s\in S'(\bar{\Q})$, letting $L:=K(s)$ we get an induced morphism $\tilde{s}:\spec(\CO_L)\rightarrow \mathfrak{S}\times_{\spec(\CO_K)}\spec{\CO_L}$, which will be a section of $\mathfrak{S}\times_{\spec(\CO_K)}\spec{\CO_L}\rightarrow \spec(\CO_L)$. Similarly for any root $\xi$ of $x$ we get sections $\tilde{\xi}:\spec(\CO_L)\rightarrow \mathfrak{S}\times_{\spec(\CO_K)}\spec{\CO_L}$. We note that the main technical feature we will need from our G-functions is the following assumption, following the discussion in \cite{andre1989g}, Ch. $X$, $\S$ $3.1$:
\begin{assumption}\label{assumprionreduction}
	Let $s\in S'(\bar{\Q})$ such that $|x(s)|_v< R_v(\mathcal{Y})$ for some finite place $v\in \Sigma_{K(s),f}$. Then $\tilde{s}$ has the same image in $\mathfrak{S}(\kappa(v))$ as $\tilde{\xi}$ for some zero $\xi$ of $x$, where $\kappa(v)$ is the residue field of $K(s)$ at $v$.
\end{assumption}

We finally record the following: 
\begin{definition}\label{vadicprox}
	Let $s\in S(L)$, with $L/K$ finite, and let $v\in \Sigma_{L}$.
	
	We say that the point $s$ is $v$-adically close to $0$ if $|x(s)|_v<\min\{1, R_v(\mathcal{Y})\}$. We furthermore say that $s$ is $v$-adically close to $\xi_t$ if furthermore $s$ is contained in the connected component of the preimage $x^{-1}(\Delta_{r_v(\mathcal{Y})})\subset S'^{an}$ that contains $\xi_t$, where $\Delta_{r_v(\mathcal{Y})}$ is the open disc, either in the rigid analytic or complex analytic sense, of radius $r_v(\mathcal{Y}):=\min\{1, R_v(\mathcal{Y})\}$.\end{definition}

	
	\section{Determining the trivial relations}\label{section:trivial}

Here we determine the so called ``trivial relations'' among the family $\mathcal{Y}$ of G-functions associated to $f'$ as in \Cref{definassgfunsglob} under \Cref{hodgegenericity}.

Throughout this section we fix a semiabelian scheme $f':\CX'\rightarrow S'$ defined over $\bar{\Q}$ satisfying the conditions of \Cref{htboundred}. In particular, we get an affine open $S:=S'\backslash \{\xi_1\ld\xi_l\}$, with $\xi_t$ the zeroes of $x$ as usual, such that $\CX:=\CX'|_S=\CE_1\times\ldots\times\CE_n$ is a product of elliptic curves over $S$. We work under the assumption that \Cref{assummonodromy}, with respect to each $\xi_t$, holds for our semiabelian scheme. It is then easy to see that for the semiabelian schemes $f'_{\lambda}:\CX'_{\lambda}\rightarrow S_1=S\cup \{\xi_1\}$ defined in \Cref{section:gcsetting}, \Cref{assummonodromy} will hold around $\xi_1$. As in \Cref{section:setting} we let $\CX_{\lambda}:=(\CX'_{\lambda})|_{S}$.

Now let us consider the semiabelian scheme $f'_{\Lambda}:\times_{\lambda\in \Lambda}\CX'_{\lambda}\rightarrow S_1$ that is the product of the semiabelian schemes $f'_{\lambda}$ above, and note that \Cref{assummonodromy} will hold around $s_1$. We furthermore fix for each $\lambda$ a basis $\{\omega^{(\lambda)}_1\ld\omega^{(\lambda)}_{2n}\}$ of $H^1_{DR}(\CX_{\lambda}/S)$, where $\omega^{(\lambda)}_{2k-1}$, $\omega^{(\lambda)}_{2k}$ are given by \Cref{gfunsthm} for the smooth coordinates  $k$ of $S_1$ and by \Cref{lemmasymplectic} for the singular coordinates $k$ of $S_1$ respectively. 

From now on let us fix an embedding $\iota:\bar{\Q}\hookrightarrow \C$. We will denote for each $\lambda\in \Lambda$ by $\Gamma^{(\lambda)}$ the basis of $R^1(f_{\lambda,\iota})_{*}(\C)$ constructed in \Cref{gfunsthm}, and write $\Gamma$ for the union of the $\Gamma^{(\lambda)}$. We then know that the relative period matrices $\CP_{\lambda}$ corresponding to each of the $f'_{\lambda}$, as well as the period matrix $\CP_{\Gamma}$ corresponding to $f'_{\Lambda}$, in a neighborhood close to $\xi_1$ in $S_\iota$ will be block diagonal with diagonal blocks given by \Cref{gfunsthm} for the smooth coordinates and by \Cref{lemmasymplectic} for the singular ones.

\subsection{Notation-Background}

We follow the general notation and ideas set out in $\S 7$ of \cite{papasbigboi} as well as $\S 6.G$ of \cite{daworr5}.

We let $m_k=1$ if $k$ is a singular coordinate for $S_1$ and $m_k=2$ if $k$ is a smooth such coordinate. We set \begin{center}
	$\B:=\A_{\bar{\Q}}^{2m_1}\times\ldots\times \A_{\bar{\Q}}^{2m_n}$, and 
\end{center}$\B_{\Lambda}=\B^{|\Lambda|}$, with $|\Lambda|$-copies of $\B$, one for each $\lambda$. We furthermore write $\spec(\bar{\Q}[X^{(\lambda)}_{i,j,k}:1\leq i,j\leq 2])=\A_{\bar{\Q}}^{2m_k}$, with $\lambda\in \Lambda$ varying, when $k$ is a smooth coordinate and $\spec(\bar{\Q}[X^{(\lambda)}_{i,1,k}:1\leq i\leq 2])=\A_{\bar{\Q}}^{2m_k}$ when $k$ is singular instead. In what follows, we alternate without mention between viewing points in these copies $\A_{\bar{\Q}}^{2m_k}$ for smooth coordinates $k$ as either $2\times 2$ matrices or just points in affine space.

Similarly we consider $\B_0:=\A_{\bar{\Q}}^4\times\ldots\times\A_{\bar{\Q}}^4$, $n$ copies, which we think of alternatively as $M_{2\times 2,\bar{\Q}}^{n}$. We let $\spec(\bar{\Q}[X^{(\lambda)}_{i,j,k}:1\leq i,j\leq 2])=\A_{\bar{\Q}}^{4}$ for each of the copies so we get a natural morphism $\B_0\rightarrow \B$ which, on the level of points, is nothing but the morphism that sends $\begin{pmatrix}
	a&b\\c&d
\end{pmatrix}\mapsto \begin{pmatrix}
a\\c
\end{pmatrix}$ for the singular coordinates and coincides with the identity for the smooth ones. We also let $\B_{0,\Lambda}$ for the product of $|\Lambda|$-copies of $\B_0$, with the obvious analogous morphism $\B_{0,\Lambda}\rightarrow \B_{\Lambda}$.

We let $P$, respectively $P_{\lambda}$, be the matrix obtained by deleting from $\CP_{\Gamma}$, respectively from $\CP_{\lambda}$, all of the columns that correspond to the $\gamma^{(\lambda)}_{2k,\iota}$ with $k$ a singular coordinate. This matrix will naturally correspond to a point in $\B_{\Lambda}( \CO_{S_\iota}(V))$, respectively in $\B( \CO_{S_\iota}(V))$, where here $V$ is a small enough analytic subset of $S_{\iota}$ as in \Cref{section:setting}. Equivalently, we may and will consider $P$, respectively $P_{\lambda}$, as a function $P:V\rightarrow \B_{\Lambda}(\C)$, respectively $P:V\rightarrow \B(\C)$.

Similarly, writing \begin{center}
	$\mathcal{Y}_{\lambda}:=\{y^{(\lambda)}_{i,j,k}:1\leq i,j\leq2, \text { }k \text{ smooth for }S_1\}\cup$
	
	$ \{ y^{(\lambda)}_{i,1,k}:i=1,2,\text{ and } k\text{ singular for }S_1\}$
\end{center}we get a corresponding point $Y_\lambda\in\B(\bar{\Q}[[x]]))$, and similarly ${Y}\in \B_{\Lambda}(\bar{\Q}[[x]]))$ corresponding to $\mathcal{Y}$. In this section our goal is to determine the equations defining the subvariety ${Y}^{\bar{\Q}[x]-Zar}$ of $(\B_{\Lambda})_{\bar{\Q}[x]}$.

Alternatively, to this family $\mathcal{Y}$ and the fixed embedding $\iota$  we can also associate a function ${Y}:V\rightarrow \B_{\Lambda}(\C)$, and similarly for the $Y_{\lambda}$. Note that for each of the smooth coordinates we will then have, from \Cref{gfunsthm}, that for all $s\in V$ \begin{equation}\label{eq:equationtrivialsmooth}
\pi_k(P_{\lambda}(s))=	\pi_k(Y_{\lambda}(s))\cdot \Pi_{\lambda,k,\iota},\end{equation}
where $\pi_k:\B\rightarrow \A^4$ denotes the projection to the $k$-th copy of $\A^4$ in $\B$.

    \subsection{The trivial relations}

We start with the following lemma, which is an analogue of Corollary $5.9$ of \cite{daworr4}.
\begin{lemma}\label{lemmagraph}
	The graph $Z'$ of the function $P:V\rightarrow \B_{\Lambda}(\C)$ is such that its $\C$-Zariski closure $Z'^{\C-Zar}\subset S_{\C}\times (\B_{\Lambda})_{\C}$ is equal to $S_{\C}\times \Theta_{1,\C}$, where $\Theta_{1,\C}$ is the subvariety of $(\B_{\Lambda})_{\C}$ cut out by the ideal\begin{equation}\label{eq:idealoverc}
		I_{0}:=\langle X^{(\lambda)}_{1,1,k}X^{(\lambda)}_{2,2,k}-X^{(\lambda)}_{1,2,k}X^{(\lambda)}_{2,1,k}-\frac{1}{2\pi i}:k\text{ is smooth for } S_1, \lambda\in \Lambda\rangle.
	\end{equation}
	
\end{lemma}
\begin{proof}We note that from the same proof as that of Lemma $6.11$ of \cite{daworr5} one can describe explicitly the $\C$-Zariski closure of the graph $Z\subset V\times \B_{0,\C}$ of the function $\CP_{\Gamma}:V\rightarrow \B_{0,\Lambda}(\C)$. Indeed, one has  that $Z^{\C-Zar}$ is equal to $S_{\C}\times \Theta_{0,\C}$, where $\Theta_{0,\C}$ is the subvariety of $\B_{0,\Lambda,\C}$ cut out by the ideal\begin{equation}\label{eq:idealovercfull}
		I_1:=\langle X^{(\lambda)}_{1,1,k}X^{(\lambda)}_{2,2,k}-X^{(\lambda)}_{1,2,k}X^{(\lambda)}_{2,1,k}-\frac{e'_k}{2\pi i}:1\leq k\leq n, \lambda\in \Lambda\rangle,
	\end{equation}where $e'_k=1$ for smooth coordinates and $e'_k=e_k$ as in part $2$ of \Cref{lemmasymplectic} for the singular coordinates.
	
	The lemma follows via the same argument as in \cite{daworr5} used to deduce their Corollary $6.12$ from their Lemma $6.11$.
\end{proof}
\begin{lemma}\label{lemmagraph2}
Let $Z_{G}$ be the graph of the function $Y:V\rightarrow \B_{\Lambda}(\C)$ and let $Z_{G}^{\C-Zar}$ be its $\C$-Zariski closure in $S_{\C}\times (\B_{\Lambda})_{\C}$. Then $Z_{G}^{\C-Zar}=S_\C\times \Theta_{\C}$ where $\Theta_{\C}$ is the subvariety of $(\B_{\Lambda})_{\C}$ cut out by the ideal\begin{equation}\label{eq:idealovercgfuns}
	I_{0}:=\langle X^{(\lambda)}_{1,1,k}X^{(\lambda)}_{2,2,k}-X^{(\lambda)}_{1,2,k}X^{(\lambda)}_{2,1,k}-1:k\text{ is smooth for } S_1,\lambda\in \Lambda\rangle.
\end{equation}
\end{lemma}
\begin{proof}Consider the automorphism $\theta:\B_{\Lambda}\rightarrow \B_{\Lambda}$ defined on the level of points $((A^{(\lambda)}_1\ld A^{(\lambda)}_n):\lambda\in \Lambda)$ by multiplying on the right by $\Pi_{\lambda,k,\iota}^{-1}$ for each $A^{(\lambda)}_k$ such that $k$ corresponds to a smooth coordinate for our curve $S_1$. 
	
By construction, see \eqref{eq:equationtrivialsmooth}, we then have that $Y=\theta\circ P$. The result follows from \Cref{lemmagraph}.\end{proof}

\begin{theorem}\label{trivialrels} With the previous notation, under \Cref{hodgegenericity}, $Y^{\bar{\Q}[x]-Zar}$ is the subvariety of $(\B_{\Lambda})_{\bar{\Q}[x]}$ cut out by 
	 \begin{equation}\label{eq:trivialsmooth}
		I_0:=\langle \det(X^{(\lambda)}_{i,j,k})-1:1\leq k\leq n, \text{ k is smooth for } S_1, \lambda\in \Lambda\rangle.
	\end{equation}
\end{theorem}
\begin{proof}
	The proof follows trivially from \Cref{lemmagraph2} since the generators of the ideal $I_0$ are all defined over $\bar{\Q}$.
\end{proof}

  	\section{Archimedean relations at CM-points}\label{section:cm}

Let us fix from now on a semiabelian scheme  $f':\CX'\rightarrow S'$ satisfying the assumptions of \Cref{htboundred}. See also \Cref{section:gcsetting} for our notation. In this section we will consider the family of G-functions associated to the above semiabelian scheme, as in \Cref{definassgfunsglob}, and construct relations among the archimedean values of this family at CM-points $s\in S(L)$, where $S$ here denotes some affine open subset of $S'\backslash\{\xi_1\ld\xi_l\}$ as in the discussion in \Cref{section:goodcovers}.\\

\begin{prop}\label{archirelationscm}
	Let $f':\CX'\rightarrow S'$ be a $G_{AO}$-admissible\footnote{Here we follow the naming convention in \Cref{gaoadmremark}.} semiabelian scheme as above. Then for any $s\in S(\bar{\Q})$ for which $\CX_s$ is CM, there exists a homogeneous polynomial $R_{s,\infty}\in L_s[X^{(\lambda)}_{i,j,k}:\lambda\in \Lambda, 1\leq i,j,\leq 2, 1\leq k\leq n]$, where $L_s/K(s)$ is a finite extension, such that the following hold: \begin{enumerate}
		
		\item $\iota_v(R_{s,\infty}(\mathcal{Y}(x(s))))=0$ for all $v\in \Sigma_{L_s,\infty}$ for which $s$ is $v$-adically close to $0$,
		
		\item\label{item2arc} $[L_s: \Q]\leq c_1(n)[K(s):K]$, where $c_1(n)$ is a constant depending only on $n$,
		
		\item $\deg(R_{s,\infty})\leq 2[L_s:\Q]$, and 
		
		\item $R_{s,\infty}(\mathcal{Y}(x))=0$ does not hold generically, in other words the relation defined by the polynomial is ``non-trivial''.	
	\end{enumerate}	
\end{prop}
\begin{definition}
	We call the field $L_s$ associated to the point $s$ the\textbf{ field of coefficients of the point }$s$. 
\end{definition}
\begin{proof}We break the proof in parts. First we create what we call ``local factors'' $R_{s,v}$, each one associated to a fixed place $v\in \Sigma_{L_s,\infty}$ for which $s$ is $v$-adically close to $0$. To do this we break the exposition into cases. First, we work under the assumption that the toric rank of the semiabelian variety $\CX'_{\xi}$, for any zero $\xi$ of $x$, is $t\geq 2$, or in other words the case where there are at least $2$ singular coordinates for our curve $S_1=S\cup\{\xi_1\}$. In the second case we will work under the assumption that there is at least one smooth coordinate which is CM, noting that in any case a singular coordinate will exist by our conventions in \Cref{section:goodcovers} as well as \Cref{defgadmi}. After this we define the polynomials $R_{s,\infty}$ in question and establish their main properties outlined in the lemma.\\
	
Before that, we fix some notation that persists in both cases. Throughout this proof we fix a point $s\in S(\bar{\Q})$ such that the fiber $\CX_s$ is CM and let $K(s)$ be its field of definition. We let $L_s$ be the compositum of the following fields \begin{enumerate}
		\item the finite extension $\hat{K}(s)/K(s)$ such that $\End_{\bar{\Q}}(\CX_s)=\End_{\hat{K}(s)}(\CX_s)$, 
		\item the CM fields $F_{k,s}:=\End_{\bar{\Q}}(\CE_{k,s})$.
	\end{enumerate}

We note that it is classically known\footnote{See the Remark on page $136$ of \cite{langelliptic}. In short, for each $1\leq i\leq n$ we have that all elements of $\End_{\bar{\Q}}(\CE_{i,s})$ will be defined over $F_{i,s}K(s)$.} here that \begin{equation}\label{eq:silverberg}
[L_s:\Q]\leq 2^n [K(s):\Q], 
\end{equation}so $c_1(n)=2^n$ in \ref{item2arc} above. Finally, let us fix a place $v\in \Sigma_{L_s,\infty}$ and let $\iota_v:L_s\hookrightarrow \C$ be the corresponding embedding. We assume from now on that $s$ is $v$-adically close to $\xi_t$ for some $1\leq t\leq l$, see \Cref{section:proximity} for the notation here.\\

\textbf{Step $0$: The issue of $v$-adic proximity} Here we follow the strategy set out in $\S6.F$ of \cite{daworr5}. Let us write $x^{-1}(\Delta_{r_v(\mathcal{Y})})=U_{v,1}\sqcup \ldots \sqcup U_{v,l}$, where $U_{v,t}$ denotes the connected component with $\xi_t\in U_{v,t}$. From the fact that $x\circ \sigma_t =x$ we have that $\sigma_t^{-1}(U_{v,t})\subset U_{v,1}$. In particular we will have that $s_t:=\sigma_t^{-1}(s)\in U_{v,1}$.

Now let $\lambda\in \Lambda$ be the minimal representative of the equivalence class of $t$, as in \Cref{section:gcsetting}, so that $\CX'_{\lambda}$ and $\CX'_t$ are generically isogenous. By Lemma $5.4$ of \cite{daworr4} we then get an isogeny $\CX_{t,s_t}\sim \CX_{\lambda,s_t}$. Since by construction we have that $\CX_{t,s_t}\simeq \CX_{s}$ we get that $\CX_{\lambda,s_t}$ is also a product of CM elliptic curves.

For ease of notation, we let from now on $f:\CG':=\CX'_\lambda\rightarrow S_1=S\cup\{\xi_1\}$ and $\CG:=\CG'|_{S}$. By abuse of notation, we will write $\CG'=\CE'_1\times\ldots \times \CE_n'$, where $\CE'_j$ are $1$-dimensional semiabelian schemes so that $\CE'_j|_{S}$ is an elliptic curve over $S$, as per our usual notation. Also for ease of notation we set $s$ for $s_t=\sigma^{-1}_{t}(s)$ for now.

As in the proof of \Cref{lemmacmperiods}, for all $1\leq k\leq n$ there exists a symplectic basis $w_{2k-1,s}$, $w_{2k,s}$ of $H^1_{DR}(\CE_{k,s}/K(s))\otimes L_s$ and a symplectic basis $\gamma'_{2k-1,s}$, $\gamma'_{2k,s}$ of $H_1(\CE_{k,s,\iota_v},\Q)\otimes L_s$ such that \eqref{eq:splitderham} and \eqref{eq:splitbetti} hold. Let us also consider the fixed basis $\{\omega_i:1\leq i\leq 2n\}$ of $H^1_{DR}(\CG/S)|_U$ and the fixed frame $\{\gamma_{j,\iota_v}:1\leq j\leq 2n\}$ of $R^1(f_{\iota_v})_{*}(\C)(1)$ chosen by the combination of \Cref{gfunsthm} and \Cref{lemmasymplectic} for our semiabelian scheme $\CG'$.

We then obtain change of bases matrices $B_{k,dR}:=\begin{pmatrix}
	a_{k,s}&b_{k,s}\\c_{k,s}&d_{k,s}
\end{pmatrix}\in\GL_{2}(L_s)$ between the bases $w_i$ and $\omega_{i,s}$ of $H^1_{DR}(\CG_{s}/L_s)$ and $B_{k,b}:=\begin{pmatrix}
\alpha_{k,s}&\beta_{k,s}\\\gamma_{k,s}&\delta_{k,s}
\end{pmatrix}\in\SL_{2}(L_s)$ between the bases $\gamma'_{j,s}$ and $\gamma_{j,s}$ of $R^1(f_{\iota_v})_{*}(\C)(1)$. Note that the fact that the entries of $B_{k,b}$ are in $L_s$ follows by construction.

Let $P_{w,\gamma',k,s}$ be the full period matrix of $\CE_{k,s}$ with respect to the bases $w_i$ and $\gamma'_j$. On the one hand we then have that \begin{equation}\label{eq:cmprep1}
	P_{w,\gamma',k,s}=B_{k,dR}\cdot \CP_{k}(s)\cdot B_{k,b},
\end{equation}where $\CP_k(s)$ denotes the value at $s$ of the relative period matrix associated to the semiabelian scheme $f'_{k,t}:\CE'_{k,t}\rightarrow C'_t$, the basis $\omega_i$, and the trivializing frame $\gamma_j$ above. On the other hand, by the construction in \Cref{lemmacmperiods} we know that \begin{equation}
	P_{w,\gamma',k,s}=	\begin{pmatrix}\label{eq:cmprep2}
		\frac{\varpi_{s,k}}{2\pi i}&0\\0&\varpi_{s,k}^{-1}
	\end{pmatrix},
\end{equation}for some transcendental number $\varpi_{s,k}$ that depends on the embedding $\iota_v$ chosen.\\
    
\textbf{Step 1: Defining the local factors} Here we break the exposition into two cases, depending on the fiber over $\xi_t$ of our semiabelian scheme.

\textbf{Case 1:} Let us assume from now on that there exist at least two singular coordinates for $S'$ and without loss of generality we assume that these are the first two. 

Let us write $B_{k,dR}\cdot \CP_{k}(s)=(p_{i,j,k})$ for convenience. We note that from our various conventions in \Cref{section:background} we know that the first column of $\CP_{k}(s)$ is actually of the form \begin{equation}
	\begin{pmatrix}
		\iota_v(y^{(\lambda)}_{1,1,k}(x(s)))\\
	    \iota_v(y^{(\lambda)}_{2,1,k}(x(s)))
	\end{pmatrix},
\end{equation}where $y^{(\lambda)}_{i,j,k}$ are members of the subfamily $\mathcal{Y}_{\lambda}$ of the family of G-functions $\mathcal{Y}$ associated to the point $\xi_1$ as in \Cref{definassgfunsglob}.

From \eqref{eq:cmprep1} and \eqref{eq:cmprep2} we get for $k=1$, $2$:\begin{equation}\label{eq:singprep1}
\begin{pmatrix}
	\frac{\varpi_{s,k}}{2\pi i}&0\\0&\varpi_{s,k}^{-1}
\end{pmatrix}=\begin{pmatrix}
\alpha_{k,s} p_{1,1,k}+\gamma_{k,s} p_{1,2,k}&        \beta_{k,s} p_{1,1,k}+\delta_{k,s} p_{1,2,k}\\ 
\alpha_{k,s} p_{2,1,k}+\gamma_{k,s} p_{2,2,k}&         \beta_{k,s} p_{2,1,k}+\delta_{k,s} p_{2,2,k}\\
\end{pmatrix}
\end{equation}

Comparing the off-diagonal elements in the equality \eqref{eq:singprep1} we get that for $k=1$, $2$ \begin{equation}\label{eq:singprep2}
  \beta_{k,s} p_{1,1,k}+\delta_{k,s} p_{1,2,k}=0\text{, and }\alpha_{k,s} p_{2,1,k}+\gamma_{k,s} p_{2,2,k}=0.
\end{equation}If for either $k=1$ or $2$ we have that $\gamma_{k,s}=0$ or $\delta_{k,s}=0$ then, since the matrix $B_{k,b}$ is invertible, we must have that $p_{1,1,k}=0$ or $p_{2,1,k}=0$. 

But by definition we have $p_{1,1,k}=\iota_v(a_{k,s}y^{(\lambda)}_{1,1,k}(x(s)) +b_{k,s}y^{(\lambda)}_{2,1,k}(x(s)))$ and $p_{2,1,k}=\iota_v(c_{k,s}y^{(\lambda)}_{1,1,k}(x(s)) +d_{k,s}y^{(\lambda)}_{2,1,k}(x(s)))$. Therefore, if for either $k=1$ or $2$, $\gamma_{k,s}=0$ or $\delta_{k,s}=0$ holds we set \begin{equation}\label{eq:singpoly1}
	R_{s,v}:=a_{k,s}X^{(\lambda)}_{1,1,k} +b_{k,s}X^{(\lambda)}_{2,1,k}\text{, or respectively }c_{k,s}X^{(\lambda)}_{1,1,k}+d_{k,s}X^{(\lambda)}_{2,1,k}.
\end{equation}

From now on let us assume that $\gamma_{k,s}$, $\delta_{k,s}\neq 0$ for $k=1$, $2$. Then \eqref{eq:singprep2} gives \begin{equation}\label{eq:singprep3}
	p_{1,2,k}=-\frac{\beta_{k,s}}{\delta_{k,s}} p_{1,1,k}\text{ and } p_{2,2,k}=-\frac{\alpha_{k,s}}{\gamma_{k,s}}p_{2,1,k}.
\end{equation}Comparing the diagonal elements in \eqref{eq:singprep1} and using \eqref{eq:singprep3} we get
\begin{equation}\label{eq:singprep4}
	\frac{\alpha_{k,s}\delta_{k,s}-\beta_{k,s}\gamma_{k,s}}{\delta_{k,s}} p_{1,1,k}=\frac{\varpi_{s,k}}{2\pi i} \text{ and }
\end{equation}\begin{equation}\label{eq:singprep5}
-	\frac{\alpha_{k,s}\delta_{k,s}-\beta_{k,s}\gamma_{k,s}}{\gamma_{k,s}} p_{2,1,k}=\varpi^{-1}_{s,k}
\end{equation}From these, together with the fact that $B_{k,b}\in \SL_2(L_s)$, we conclude that for $k=1$, $2$\begin{equation}\label{eq:singprep6}
p_{1,1,k}\cdot p_{2,1,k}= -\gamma_{k,s}\delta_{k,s}\frac{1}{2\pi i}.
\end{equation}

Finally, from \eqref{eq:singprep6} we can get rid of the $2\pi i$ to conclude that \begin{equation}
	\gamma_{2,s}\delta_{2,s}p_{1,1,1}\cdot p_{2,1,1} =\gamma_{1,s}\delta_{1,s} p_{1,1,2}\cdot p_{2,1,2}.
\end{equation}As we have seen above, we can then associate to the place $v$ and the point $s$ the polynomial 
\begin{equation}\label{eq:singpoly2}
	\begin{split}
R_{s,v}:=\gamma_{2,s}\delta_{2,s}(a_{1,s}X^{(\lambda)}_{1,1,1} +b_{1,s}X^{(\lambda)}_{2,1,1})(c_{1,s}X^{(\lambda)}_{1,1,1}+d_{1,s}X^{(\lambda)}_{2,1,1})\\
-\gamma_{1,s}\delta_{1,s}(a_{2,s}X^{(\lambda)}_{1,1,2} +b_{2,s}X^{(\lambda)}_{2,1,2})(c_{2,s}X^{(\lambda)}_{1,1,2}+d_{2,s}X^{(\lambda)}_{2,1,2}).\end{split}\end{equation}

We note that in either case $R_{s,v}$ is homogeneous of degree at most $2$ and that $\iota_v(R_v(\mathcal{Y}(x(s)))=0$.\\
    
\textbf{Case 2:} Let us now assume that there exists at least one smooth coordinate for $S'$ that is CM and without loss of generality assume that it is the first one. 

Again combining \eqref{eq:cmprep1} and \eqref{eq:cmprep2} for $k=1$, together with the description of $\CP_{k,v}(s)$ given by \Cref{gfunsthm}, we conclude that 
 \begin{equation}\label{eq:cmprep3}
	\begin{pmatrix}
		\frac{\varpi_{s,1}}{2\pi i}&0\\0&\varpi_{s,1}^{-1}
	\end{pmatrix}=\iota_v\left(  \begin{pmatrix}
		a_{1,s}&b_{1,s}\\c_{1,s}&d_{1,s}
	\end{pmatrix} Y_{G,k}(x(s))  \begin{pmatrix}
		\frac{\varpi_{0,1}}{2\pi i}&0\\0& \varpi_{0,1}^{-1}
	\end{pmatrix} \begin{pmatrix}
		\alpha_{1,s}&\beta_{1,s}\\ \gamma_{1,s}&\delta_{1,s}
	\end{pmatrix}\right),
\end{equation}noting that $\varpi_{s,1}$ itself depends on the embedding $\iota_v$.

As before for convenience let us write $(p_{i,j}):= B_{1,dR}\cdot Y_{G,1}(x(s))$. Rewriting \eqref{eq:cmprep3} we get 
\begin{equation}\label{eq:cmprep4}
	\begin{pmatrix}
		\frac{\varpi_{s,1}}{2\pi i}&0\\0&\varpi_{s,1}^{-1}
	\end{pmatrix}=\iota_v\left(\begin{pmatrix}
	p_{1,1}\alpha_{1,s}\frac{\varpi_{0,1}}{2\pi i}+p_{1,2}\gamma_{1,s}  \varpi_{0,1}^{-1}&	p_{1,1}\beta_{1,s}\frac{\varpi_{0,1}}{2\pi i}+p_{1,2}\delta_{1,       s}  \varpi_{0,1}^{-1}\\
		p_{2,1}\alpha_{1,s}\frac{\varpi_{0,1}}{2\pi i}+p_{2,2}\gamma_{1,s}  \varpi_{0,1}^{-1}&	p_{2,1}\beta_{1,s}\frac{\varpi_{0,1}}{2\pi i}+p_{2,2}\delta_{1,s}  \varpi_{0,1}^{-1}
\end{pmatrix}\right).
	\end{equation}

Considering the equalities given from the off-diagonal entries in \eqref{eq:cmprep3} we conclude that \begin{equation}
	A\frac{\varpi_{0,1}}{2\pi i}+B\varpi_{0,1}^{-1}=0\text{ and } C\frac{\varpi_{0,1}}{2\pi i}+D\varpi_{0,1}^{-1}=0,
\end{equation}where $\begin{pmatrix}
A&B\\C&D
\end{pmatrix}= \begin{pmatrix}p_{1,1}\beta_{1,s}&p_{1,2}\delta_{1,s}\\p_{2,1}\alpha_{1,s}&p_{2,2}\gamma_{1,s}\end{pmatrix}$. From this we get that $\det\begin{pmatrix}
A&B\\C&D\end{pmatrix}=0$. In this equation, we replace $p_{i,j}$ by its expression in terms of the entries of $B_{1,dR}$ and $Y_{G,1}(x(s))$. Using the fact\footnote{We note here that in the case of smooth coordinates all our chosen bases in the de Rham side may be chosen to be symplectic so that $B_{1,dR}\in\SL_2(L_s)$ in this case. This allows us to simplify our computations without having to keep track of $\det B_{1,dR}$.} that $\det B_{1,dR}=\det B_{1,b}=1$, we get the relation\begin{equation}\label{eq:cmpoly1}
	\begin{split} \iota_v(a_{1,s}c_{1,s}y^{(\lambda)}_{1,1,1}(x(s))y^{(\lambda)}_{1,2,1}(x(s))+b_{1,s}d_{1,s}y^{(\lambda)}_{2,1,1}(x(s))y^{(\lambda)}_{2,2,1}(x(s))\\
	+(2b_{1,s}c_{1,s}+1)y^{(\lambda)}_{1,1,1}(x(s))y^{(\lambda)}_{2,2,1}(x(s))\\
	-(1+b_{1,s}c_{1,s}+\beta_{1,s}\gamma_{1,s})\det(y^{(\lambda)}_{i,j,1}(x(s))))=0.
	\end{split}
\end{equation}

This will naturally correspond to a polynomial $R_{s,v}\in \bar{\Q}[X^{(\lambda)}_{i,j,k}]$ as in the previous case. Note that by construction we will have that $\iota_v(R_{s,v}(\mathcal{Y}(x(s))))=0$. Note also that $R_{s,v}$ is homogeneous of degree $2$. This last fact is easy to see once one writes $R_{s,v}$ as a sum of monomials, upon which step the fact that $\det B_{1,dR}=\det B_{1,b}=1$ makes it impossible that all coefficients of the polynomial in question are zero.\\
    
\textbf{Step 3: The polynomial} $R_{s,\infty}$ Let us now consider the following polynomial\begin{equation}\label{eq:finalpolyy1cm}
	R_{s,\infty}(X^{(\lambda)}_{i,j,k}):=\prod_{\underset{s \text{ is } v\text{-adically close to }0}{v\in \Sigma_{L_s,\infty}}}^{}R_{\sigma_{t(v)}^{-1}(s),v}(X^{(\lambda)}_{i,j,k}),
\end{equation}where\footnote{Here we note that in each case the $t\in\{1\ld l\}$ for which $s$ is $v$-adically close to $\xi_t$, will depend on the place $v$, hence the notation used here.} $R_{\sigma_{t(v)}^{-1}(s),v}(X^{(\lambda)}_{i,j,k})$ are the polynomials in \eqref{eq:singpoly1} or \eqref{eq:singpoly2}, depending on the cases outlined in the first case we examined, or the polynomials corresponding to \eqref{eq:cmpoly1}.

We note that by construction we will have that $\deg R_{s,\infty}\leq 2[L_s:\Q]$ and hence statement $(3)$ of the Lemma follows. We also note that by construction of the local factors $R_{s,v}$ statement $(1)$ of our Lemma holds as well. Indeed, let $v\in \Sigma_{L_s,\infty}$ be such that $s$ is $v$-adically close to $0$. Indeed, letting $t=t(v)$ in the above notation and $s_t=\sigma_t^{-1}(s)$ as earlier, since $x\circ \sigma_t^{-1}=x$ we get that $y^{(\lambda)}_{i,j,k}(x(s))=y^{(\lambda)}_{i,j,k}(x(s_t))$. But then by construction we know that $\iota_v(R_{s_t,v}(y^{(\lambda)}_{i,j,k}(x(s_t)))))=0$.\\

\textbf{Step 4: Non-triviality} The only thing we are left with showing is statement $(4)$ of the Lemma. This would show the ``non-triviality'' of the relation among the values at $x(s)$ of the G-functions of our family $\mathcal{Y}$ in the notation of \cite{andre1989g} Ch. $VII$, $\S$ $5$. In practical terms, we need to show that $R_{s,\infty}\notin I_{0}$, where $I_0$ is the ideal described in \Cref{trivialrels}.

We note here that the ideal $I_0$ is clearly prime. This can be seen for example by the fact that the subvariety $V(I_0)\subset (\mathbb{B}_{\Lambda})_{\bar{\Q}[x]}$ it defines is irreducible. Indeed, it is easy to see that $V(I_0)$ is isomorphic to $\mathbb{A}^{2|\Lambda| t}\times \SL_2^{(n-t)|\Lambda|}$, where $t$ here denotes the toric rank of any of the fibers $\CX'_{\xi_j}$. It thus suffices to show that none of the local factors $R_{s,v}$ can be in $I_0$. 

Note that the local factors are such that only the $X^{(\lambda)}_{i,j,k}$ for a single $\lambda\in \Lambda$ will appear in $R_{s,v}$. If we had $R_{s,v}\in I_0$ for some $v$, writing $\lambda\in\Lambda$ for the unique element above, upon quotienting out by the ideal
\begin{center}
	$I_0':=\langle X^{(\lambda')}_{i,1,k}, X^{(\lambda')}_{i,i,k'}-1, X^{(\lambda')}_{i,j,k'};\lambda'\neq\lambda,i\neq j, k' \text{ smooth, } k\text{ singular}\rangle$,
\end{center}we would have that $R_{s,v}\in I_\lambda$, where $I_{\lambda}:=\langle \det(X^{(\lambda)}_{i,j,k})-1:1\leq k\leq n, \text{ k is smooth for } S_1\rangle$. Thus we write $X_{i,j,k}$ for the $X_{i,j,k}^{(\lambda)}$, dropping any reference to $\lambda$, to simplify our notation for the remainder of this proof. 

First let us assume that $R_{s,v}$ is of the form \eqref{eq:singpoly1}. It is trivially seen that $R_{s,v}\neq 0$ since $B_{k,dR}\in \SL_2(L_s)$. Assume without loss of generality that $R_{s,v}=a_{1,s}X_{1,1,1} +b_{1,s}X_{2,1,1}$ with $a_{1,s}\neq 0$. Then it is trivial to see that we cannot have $R_{s,v}\in I_0$ since $I_0$ is generated by the polynomials $g_k:=\det(X_{i,j,k})-1$ where $1\leq k\leq n$ runs through the smooth coordinates for $S'$, and in this case $k=1$ is a singular coordinate.

Now let us assume that $R_{s,v}$ is as in \eqref{eq:singpoly2}, without loss of generality assuming that the two singular coordinates are $k=1$ and $k=2$. Then we have $\gamma_{k,s}\neq 0$ and $\delta_{k,s}\neq 0$ for $k=1$, $2$ by assumption in this case and again the fact that $B_{k,dR}\in \SL_{2}(L_s)$ shows that $R_{s,v}\neq 0$. It is easy to see once again by the above argument that we cannot have $R_{s,v}\in I_0$. 

Finally, let us assume that we are in the case where $R_{s,v}$ is the polynomial that corresponds to \eqref{eq:cmpoly1}, without loss of generality assuming that $k=1$ is a CM coordinate for $S'$. Let $g_{k}=\det(X_{i,j,k})-1$ with $k$ running through the smooth coordinates of $S'$ as usual. Then, it is easy to see that $R_{s,v}\in I_\lambda$ implies that $R_{s,v}\in (g_1)\leq L_s[X_{i,j,1}:1\leq i,j,\leq 2]$. Since $(g_1)\subset m_{1}:=\langle  X_{1,1,1}-1, X_{1,2,1},X_{2,1,1},X_{2,2,1}-1  \rangle $ we must have $R_{s,v}\in m_1$ which is easily seen to imply $R_{s,v}\left(\begin{pmatrix}
	1&0\\0&1
\end{pmatrix}\right)=2b_{1,s}c_{1,s}+1=0$. 

On the other hand letting $m_N:=\langle X_{1,1,1}-N, X_{1,2,1}-1,X_{2,1,1}+\frac{1}{2},X_{2,2,2}-\frac{1}{2N}\rangle$ for all $N\in \N$, $N\geq 2$, and noting that $(g_1)\subset m_N$, we will have that $R_{s,v}\in m_N$ for all $N\geq 2$, $N\in \N$. Keeping in mind that $2b_{1,s}c_{1,s}+1=0$ we get that \begin{equation}
4	a_{1,s}c_{1,s}N^2-b_{1,s}d_{1,s}=0
\end{equation}for all $N$ as above. This gives $a_{1,s}c_{1,s}=b_{1,s}d_{1,s}=0$ which, together with $2b_{1,s}c_{1,s}+1=0$, is impossible since $B_{1,dR}\in\SL_2(L_s)$.\end{proof}

    
    \section{Proof of the height bounds}\label{section:heightbounds}

We use arguments centered around Gabber's lemma, as in \cite{andre1989g} and \cite{papaszp}, to rule out p-adic proximity of the points we are interested in to the point $s_0$. After this we finally come to the proof of the height bounds we want.\\

\subsection{$p$-adic proximity}

\begin{lemma}\label{padicproxao}Let $f':\CX'\rightarrow S'$ be a $G_{AO}$-admissible\footnote{Here we follow the convention in \Cref{gaoadmremark} again.} semiabelian scheme satisfying the assumptions of \Cref{htboundred}. Let $s\in S(\bar{\Q})$ be a CM point with field of coefficients $L_s$ defined as in \Cref{archirelationscm}. Furthermore assume that there exists at least one singular coordinate for $S'$.
	
	Then if $v\in \Sigma_{L_s,f}$ is some finite place of $L_s$, the point $s$ is not $v$-adically close to $s_0$.
\end{lemma}
\begin{proof}
	Using \Cref{assumprionreduction}, the proof of Lemma $5.4$ in \cite{papaszp} shows that if $s$ was $v$-adically close to $s_0$ then the special fiber of the connected N\'eron model of $\CX_{s}\times _{K(s)}L_{s,v}$ would be the same as that of $\CX_{0}\times _{K}L_{s,v}$. 
	
	Since each coordinate $\CE_{k,s}$ is CM it will have potentially good reduction at $v$, by Theorem $6$ of \cite{serretate}, while for $\CX_{0}$ we know that at least one of the coordinates is isomorphic to $\G_m$ which is a contradiction.
\end{proof}

\subsection{Proof of the height bound}

Following our discussion in \Cref{section:goodcovers} it suffices to establish the following:
\begin{prop}\label{propredhtbound}
	\Cref{heightboundandreoort} holds for a semiabelian scheme satisfying the assumptions in \Cref{htboundred}.
\end{prop}

\begin{proof}In the construction of the bases $\omega_i$ of $H^1_{DR}(\CX^{(\lambda)}/S)|_{U_\lambda}$ we have excluded a finite number of points, i.e. the points in $S \backslash U_\lambda$. Let $M:=\max\{ h(x(P)):P\in  (S\backslash U_\lambda )(\bar{\Q}), \lambda\in\Lambda \}$. 

Now fix a point $s$ for which $\CX_s$ is CM and let \begin{center}
	$\Sigma(s):=\{v\in\Sigma_{L_s,\infty}:s \text{ is } v\text{-adically close to }s_0  \}$. 
\end{center}If $\Sigma(s)=\emptyset$ then as in the proof of Theorem $1.3$ of \cite{papasbigboi}, see $\S$ $12$ there, we know that \begin{center}
	$h(x(s))\leq \rho(\mathcal{Y}):=\underset{1\leq t\leq l}{\max}\rho(\mathcal{Y}_{\lambda})$.
\end{center}
	
On the other hand, if $\Sigma(s)\neq \emptyset$ combining \Cref{archirelationscm} with \Cref{padicproxao} we get non-trivial and global relations among the values of our G-functions at $x(s)$, in the terminology of Ch. $VII$, $\S$ $5$ of \cite{andre1989g}. Thus, the ``Hasse principle'' of Andr\'e-Bombieri, CH. $VII$, Theorem $5.2$ in \cite{andre1989g}, gives that \begin{equation}
		h(x(s))\leq c_{0,1}\deg(R_{s,\infty})^{c_2}.
	\end{equation}We note that the constant $c_{0,1}$ will only depend on the differential operator $L_{\lambda}$ associated via the Gauss-Manin connection with our choice of bases and the family of G-functions $\mathcal{Y}$, while the constant $c_2$ will only depend on $n$ and $|\Lambda|$. Finally, since $[L_s:\Q]\leq c_1(n) [K(s):\Q]$, by \Cref{archirelationscm}, the result follows, setting $c_1:=\max(c_{0,1}c_1(n)^{c_2}, M, \rho(\mathcal{Y}))$.
\end{proof}

   
    \section{Effective Brauer-Siegel}\label{section:applications}

We are finally able to establish \Cref{maintheorem} pairing our \Cref{heightboundandreoort} together with results of Masser-W\"{u}stholz.

\begin{proof}[Proof of \Cref{maintheorem}]
	This proof is pretty much verbatim that of Proposition $5.13$ of \cite{daworr4}. Throughout let us fix a point $s$ as in the statement.
	
Let us fix a compactification $\bar{Z}$ of $Z$ in $X(1)^n\simeq (\mathbb{P}^1)^n$. Then we can find a finite \'etale cover of $\bar{Z}$, $g:\bar{S}\rightarrow  \bar{Z}$, such that after possibly base changing by a finite extension $K'/K$, we have that the semiabelian scheme $f':\CX'\rightarrow S'$, where \begin{enumerate}
	\item $S'$ is an open subset of $\bar{S}$ such that $g(S')\cap (X(1)^n\backslash Y(1)^n)=\{z_0\}$ with preimage $s_0\in S'(K)$, 
	
	\item $f:\CX=\CE_1\times\ldots \times \CE_n \rightarrow S'\backslash \{s_0\}$ is the pullback of the universal family, and 
	
	\item $f':\CX'\rightarrow S'$ is the connected N\'eron model of $f$ over $S'$,
\end{enumerate}is such that it satisfies \Cref{assummonodromy}, \Cref{assumptionfield}, and \Cref{hodgegenericity}.

We can then apply \Cref{heightboundandreoort} for any $\tilde{s}\in C_4(\bar{\Q})$ that is a preimage of $s$ in the cover $C_4$ of $S'$ described in \Cref{htboundred}. Then we know that \begin{equation}\label{eq:htboundao}
		h(x(\tilde{s}))\leq c_1[K(s):\Q]^{c_2},
\end{equation}with the constants that appear here being independent of the point $s$.

Letting $\rho_i$ be the compositions $S\xrightarrow{g}Z\xrightarrow{\pi_i }Y(1)\simeq \A_1$, and applying \cite{silverheights} Proposition $2.1$ we get that for all $1\leq k\leq n$ we have \begin{equation}\label{eq:silverhtcomp}
| h(\rho_k(s))-12 h_F(\CE_{k,s})|\leq c_3 \log \max\{2, h(\rho_k(\tilde{s}))\}.
\end{equation}Note here that the constant $c_3$ is just a constant independent of our setting.

On the other hand, we have from standard facts about Weil heights that\begin{equation}\label{eq:weilheightmach}
		|h(x(\tilde{s}))- c_5 h(\rho_k(\tilde{s}))|\leq c_6 h(x(\tilde{s})),
\end{equation}here $c_5$ and $c_6$ will depend on our curve. 

On the other hand note that from \cite{mwendoesti} we know that for all $1\leq k\leq n$ we have \begin{equation}\label{eq:masserwu}
	|\disc(\End_{\bar{\Q}}(\CE_{k,s}))|\leq c_7 \max\{[K(s):\Q],h_F(\CE_{k,s})\}^{c_8}
\end{equation}where $c_7$ and $c_8$ are positive constants that are also independent of our setting.

Combining \eqref{eq:silverhtcomp} together with \eqref{eq:weilheightmach} and \eqref{eq:masserwu} we conclude that there exist constants $c_9$, $c_{10}$ independent of our chosen point $s$ such that for all $1\leq k\leq n$ we have\begin{equation}\label{eq:aopenult}
|\disc(\End_{\bar{\Q}}(\CE_{k,s}))|\leq c_9 \max\{[K(s):\Q],h(x(\tilde{s}))  \}^{c_{10}}.
\end{equation}Pairing this last equation with \eqref{eq:htboundao} we have concluded the proof.
\end{proof}

\begin{remark}\label{remarkonconstants}We close off with some remarks on the final constants that appear in \Cref{maintheorem}.
	
	We start with the ``exponent constant'' $C_2$, which at the end of the day will be of the form $\frac{1}{c_2\cdot c_8}$, where $c_2$ is the constant in \Cref{heightboundandreoort} and $c_8$ is as in \eqref{eq:masserwu}. The constant $c_2$ that appears in \Cref{heightboundandreoort}, after \Cref{htboundred}, will be $\leq 6\cdot |\Lambda|-2$. This can be seen by the fact that the family $\mathcal{Y}$ will satisfy a differential system in the notation of Theorem $5.2$, Ch. $VII$ of \cite{andre1989g}, with $\mu=2n\cdot |\Lambda|$. Indeed this will come out of the Gauss-Manin connection of the semiabelian scheme $f'_{\Lambda}$ considered in the beginning of \Cref{section:cm}. On the other hand, the constant $c_8$ can be crudely bounded by $c_8\leq 8\cdot 768$, see for example the bottom of page $650$ of \cite{mwendoesti}, where in our case $l=2$ and $n=1$ in the notation there. Finally, we note that the quantity $|\Lambda|$ can be crudely bounded in terms of the genus of our curve $S'$, see here the original construction in Lemma $5.1$ of \cite{daworr4}.
	
	Let us now look at the constant $C_1$. Again this will be nothing but $\frac{1}{c_1c_{11}}$, where $c_1$ is as in \eqref{eq:htboundao} and $c_{11}$ is polynomially dependent on the constants $c_5$ and $c_7$ of the previous proof. In the notation of \cite{andre1989g}, and in particular $VI.4$, we let $L_\lambda:=\frac{d}{dx}-\Gamma_\lambda$ be the differential operators associated to each of the families $\mathcal{Y}_{\lambda}$ of \Cref{definassgfunsglob}. We have already given a crude description of the constant $c_1$ that appears in \eqref{eq:htboundao} at the end of the proof of \Cref{propredhtbound}. The constant $c_{0,1}$ that appears there will depend\footnote{See the Proposition on page $133$ of \cite{andre1989g} fore more details here.} on $n\cdot|\Lambda|$, $\rho(\mathcal{Y})$, $\sigma(\mathcal{Y})$, $|Sin (L_\lambda)|$, i.e. the number of singularities of $L_\lambda$, and $\sigma(L_\lambda)$ the size of $L_\lambda$. For a definition of the latter quantity see \cite{andre1989g} Ch. $IV$. By the Theorem on page $123$ of \cite{andre1989g} one can replace the dependence on $\sigma(L_{\lambda})$ by a dependence on $\sigma(\mathcal{Y}_{\lambda})$ and the quantity ``$s$'' defined on page $120$ in \cite{andre1989g} that depends on the degrees of the denominators and numerators of the entries of the matrix $\Gamma_\lambda$ associated to the bases $\omega_i$ via the Gauss-Manin connection. On the other hand, the constant $c_7$, coming from the work \cite{mwendoesti} of Masser-W\"ustholz, may be taken to be an absolute constant, denoted by ``$c_5$'' on page 650 of \cite{mwendoesti}. All in all, $C_1$ will be a constant depending on the pair $(C',s_0)$ via quantities coming from the Gauss-Manin connection, or other geometric quantities of our curve, e.g. its genus.
 \end{remark}

 \textbf{Acknowledgments:} The author thanks Yves Andr\'e for answering some questions about his work on G-functions and Chris Daw for interesting conversations around the G-functions method. The author also heartily thanks the anonymous referee for their careful reading of an earlier draft of this article, for their suggestions which greatly improved the exposition, as well as for the correction of multiple mistakes on the author's part.\\

\textbf{ Funding:} This work was funded by Michael Temkin's ERC Consolidator Grant 770922 - BirNonArchGeom. During the revision process the author was supported by the Minerva Research Foundation Member Fund while in residence at the Institute for Advanced Study for the academic year 2025-26.\\
    
\textbf{Statement of interests: }The authors have no competing interests to declare that are relevant to the content of this article.
	\bibliographystyle{alpha}
	
	\bibliography{biblio}
\end{document}